%% file: draft.tex
\documentclass{article}

\usepackage[utf8]{inputenc}
\usepackage[T1]{fontenc}
\usepackage[english]{babel}
\usepackage{enumitem}
\usepackage[hidelinks]{hyperref}
\usepackage{graphicx,caption,subcaption}
\usepackage{natbib}
\usepackage[retainorgcmds]{IEEEtrantools}
\usepackage{amsmath,mathrsfs,amsthm,amssymb,mathtools}
\usepackage[textwidth=384pt,textheight=600pt]{geometry}

\def\Id{\mathbf{I}}
\def\0{\mathbf{0}}
\def\cross[#1]{\left[ #1 \times \right]}

\input{notation}
\graphicspath{{images/}}

\usepackage{xcolor}
\usepackage{floatrow}
\usepackage{booktabs}
\usepackage{tabularx}
\newcolumntype{C}{>{\centering\arraybackslash}X}

\definecolor{gray1}{gray}{0.92}
\definecolor{gray2}{gray}{0.85}
\definecolor{gray3}{gray}{0.78}
\definecolor{gray4}{gray}{0.71}

\title{Relative Equilibria of Magnetic Micro-Swimmers}
\author{Pauline Rüegg-Reymond\footnote{Institute of Mathematics, \'Ecole polytechnique f\'ed\'erale de Lausanne (EPFL) Station 8, CH-1015 Lausanne} \and Thomas Lessinnes\footnote{\'Ecole Polytechnique de Bruxelles, ULB. Email address for correspondence: Thomas.Lessinnes@ulb.ac.be}
}

\newtheorem{longtime behaviour steadystate}{Proposition}[section]
\newtheorem{longtime behaviour dynamics}[longtime behaviour steadystate]{Proposition}

\begin{document}
\maketitle

\begin{abstract}
We revisit the dynamics of a permanent-magnetic rigid body submitted to a spatially-uniform steadily-rotating magnetic field in Stokes flow. We propose an analytical parameterisation of the full set of equilibria depending on two key experimental parameters, and show how it brings further understanding that helps to optimise magnetisation and operating parameters. The system is often bistable when it reaches its optimal swimming velocity. A handling strategy is proposed that guarantees that the correct equilibrium is reached. 
\end{abstract}

\section{Introduction}

Over the past three decades, the engineering community has imagined, built, and tested artificial micro-swimmers inspired by micro-organisms equipped with propelling helical flagella~\citep[]{Honda1996,Dreyfus2005,Bell2007,Ebbens2010}. These devices are destined to various environmental and biomedical applications~\citep{Nelson2010,Gao2014,Yan2017,Bente2018}, including water decontamination~\citep[]{Mushtaq2015}, targeted drug delivery~\citep[]{Felfoul2016}, and enhanced sperm motility~\citep[]{Medina-Sanchez2016}.

The main idea is to rotate a helical body about its helical axis, and to count on the fluid-structure interaction (in the Stokes limit) to obtain a net translation of the chiral shape along its helical axis. To achieve this, the swimmer is often magnetised in a direction $\mlab$ that is perpendicular to its helical axis~\citep[]{Ghosh2009,Tottori2012,Ghosh2012,Peyer2013b}. To obtain motion in a fixed lab direction $\elab_3$, an external magnetic field $\Blab$ is applied perpendicularly to $\elab_3$ and then rotated about $\elab_3$. If the rotation rate of the external field is slow enough, the magnetic moment $\mlab$ has time to align with $\Blab$ and to follow it as it rotates about $\elab_3$. If all goes well, the helical axis spontaneously aligns with $\elab_3$ and translation in that direction does occur.

The Mason number $\Ma$ is the key nondimensional parameter that is proportional to the angular velocity of the external field. It turns out that if $\Blab$ rotates too slowly, the swimmers tumbles instead of swimming~\citep[]{Tottori2012,Ghosh2012,Ghosh2013,Morozov2014}, namely the helical axis does not align with $\elab_3$ and the axial velocity remains small. The dependance of the wobbling angle $\beta$ between the helical axis and $\elab_3$ was studied in detail by \cite{Man2013} and \cite{Morozov2014}. As one increases $\Ma$, a bifurcation occurs to the swimming mode for which $\beta$ is small, and good axial velocities can be achieved. If $\Ma$ is increased further, a secondary bifurcation occurs, where the magnetic moment of the swimmer can no longer follow the external field, a phenomenon that is referred to as step-out \citep{Mahoney2014}. The swimming mode breaks down and the axial velocity falls sharply with further increase of $\Ma$ beyond step-out.

Although the classical setup described above provides a successful solution for driving the swimmer, it does face some practical limitations. In particular, it is difficult to precisely control the magnetisation of swimmers -- see~\cite{Morozov2014} for a discussion. Even if the magnetisation was exactly controlled, because experiments deal with finite helices exhibiting rough surfaces, the optimal magnetisation is rarely exactly perpendicular to the helical axis. Whenever the magnetisation is not perpendicular to the easy-rotation axis, it is not clear that the external field should rotate perpendicularly to the intended direction of translation. In fact, we will show that it should not. Furthermore, it is well known that non helical swimmers can also be made to swim \citep[]{Meshkati2014,Fu2015,Vach2015,Morozov2017}. It is particularly remarkable that \cite{Morozov2017} showed both that a deformed helix swims faster than an actual helix and that it is possible to make mirror-symmetric  (zero chirality) objects swim. In all those cases, there is no reason to limit $\Blab$ to being perpendicular to $\elab_3$. To our knowledge, \cite{Meshkati2014} presented the only investigation (besides the present text) that examined an arbitrary angle $\psi$ between $\Blab$ and $\elab_3$. 

The present study considers a rigid, neutrally buoyant swimmer of arbitrary shape that is a permanent magnet driven by a spatially uniform magnetic field that rotates and at constant angle $\psi$ from a fixed direction $\elab_3$. We provide the first fully analytic parameterisation of the set of relative equilibria of the system. It allows to prove that this system has almost always 0, 4 or 8 relative equilibria (in practice we find that at most two are stable). It also yields an explicit parameterisation of the achievable axial velocities of a given shape with a specified magnetisation. It then becomes particularly easy to optimise the driving parameters $\Ma$ and $\psi$ so as to achieve maximal axial velocity, which is often found for $\psi$ well away from the habitual $\pi/2$. Furthermore, precise optimal results on the magnetisation are also discussed. Finally, the swimming regime often occurs when the system is bistable. We propose three experimental strategies to guarantee that the swimmer ends up in the faster stable equilibrium.

The document is structured as follows. In section~\ref{sec-physics}, we revisit a derivation of the equations of movement. It is included here because there is a pleasant symmetry between these magnetic swimmers, which are driven by a constant torque (at relative equilibrium), and the sedimentation problem studied by \cite{Gonzalez2004}, where the object is driven by a constant force (at relative equilibrium). The resulting differential problem is summarised in equations~(\ref{cl}-\ref{outer layer system}). It consists in a nonlinear ODE on $SO(3)$ that is equivalent to the problem studied by \cite{Meshkati2014}. Section~\ref{section: analysis}, contains the analytical treatment mentioned above. Section~\ref{section: numerics} showcases the tools developed so far by studying two helical swimmers in details. Section~\ref{sec: comparison} revisits the three-beads cluster to offer a comparison with other studies and in particular with the results of \cite{Meshkati2014} and \cite{Morozov2017}. More examples are also available in the thesis of \cite{Ruegg-Reymond2019}.

\section{Modelling the Dynamics}
\label{sec-physics}

\subsection{Rigid Body Kinematics and Dynamics}

The configuration of an arbitrary rigid body is given by its position $\xlab$ and its orientation, represented by a right-handed, orthonormal and material frame $R = \begin{bmatrix} \dlab_{1} & \dlab_{2} & \dlab_{3} \end{bmatrix}$ - the body frame. The kinematic of the body is specified by $\vlab$ and $\omegalab$ which are respectively its linear and angular velocities:
\begin{align}
	\label{kinematics}
	\dot{\xlab} &= \vlab, & \dot{R} &= \cross[\omegalab] R. 
\end{align}
Note that the second equation in\re{kinematics} can also be casted in the form $\dot{\dlab}_{i} = \omegalab \times \dlab_{i}$ for $i = 1,2,3$.

The linear $\plab$ and angular $\Llab$ momenta about the centre of mass are
\begin{align}
	\label{velocity momentum}
	\plab &= k \vlab, & \Llab &= \inertlab \omegalab,
\end{align}
where $k$ is the mass of the body and $\inertlab = \inertlab^T > 0$ its inertia tensor with respect to the centre of mass; it is constant in the body frame.

The rigid body dynamic is governed by the balance of linear and angular momenta
\begin{align}
	\label{balance momentum}
	\dot{\plab} &= \flab, & \dot{\Llab} &= \taulab,
\end{align}
where $\flab$ and $\taulab$ are respectively the resultant force and torque acting upon the body. In our case, $\flab$ and $\taulab$ arise from hydrodynamic drag, applied magnetic field, gravity and buoyancy. We first focus on the role played by the drag.

We make the assumption that the body is in an unbounded fluid at low Reynolds number. This implies in particular that the loads due to hydrodynamic drag depend linearly on the velocities:
\begin{equation}
	\label{hydrodynamic loads}
	\begin{bmatrix}
		\flab\dsup \\ \taulab\dsup
	\end{bmatrix}
	= - \Drag
	\begin{bmatrix}
		\vlab \\ \omegalab
	\end{bmatrix}
	= -
	\begin{bmatrix}
		\draglab_{11} \vlab + \draglab_{12} \omegalab \\
		\draglab_{12}^T \vlab + \draglab_{22} \omegalab
	\end{bmatrix},
\end{equation}
where $\Drag \in \mathbb R^{6\times6}$ is the symmetric and positive definite drag tensor the components of which are constant in the body frame~\citep[]{Happel1983}.

We gather all forces and torques from effects other than drag respectively in  $\flab\extsup$ and $\taulab\extsup$ and expand \eqref{balance momentum} in the body frame:
\begin{equation}
	\label{balance body generic}
	\begin{aligned}
		\dot{\pbody} + \omegabody \times \pbody = &-\dragbody_{11} \vbody - \dragbody_{12} \omegabody + \fbody\extsup \\
		\dot{\Lbody} + \omegabody \times \Lbody = &-\dragbody_{12}^T \vbody - \dragbody_{22} \omegabody + \taubody\extsup
	\end{aligned}
\end{equation}
where the sans-serif letter $\pbody$ denotes the triple $\left( \plab \cdot \dlab_{i}  \right)_{i=1,2,3}$, the entries of $\dragbody_{11}$ are $\left( \dlab_{i} \cdot \draglab_{11} \dlab_{j} \right)_{i,j = 1,2,3}$, and so on. In particular, the matrices $\dragbody_{ij}$ appearing in\re{balance body generic} are constant.

We will later focus on the specific case where the loads $\fbody\extsup$ and $\taubody\extsup$ arise from gravity, buoyancy, and an interplay between a magnetic moment of the body and an external magnetic field. Note however that\re{balance body generic} is valid for arbitrary external loads and so is the dimensional analysis carried out in section~\ref{section: dim analysis}.

\subsection{Dimensional Analysis}
\label{section: dim analysis}

The relevant scales for our problem are a characteristic body length $\ell$, mass $k$ and volume $V$, the dynamic viscosity of the fluid $\eta$, its mass density $\rho_{f}$, and the magnitude of the applied force $N$. The characteristic time scale for the resulting settling phenomenon is then $t_c= \frac{\eta \, \ell^2}{N}$ (see~\cite{Gonzalez2004} for a discussion). 

Accordingly, we non-dimensionalise according to
\begin{align*}
\tbar &= \frac t {t_c}, &\vbodybar &= \frac{t_{c}}{\ell} \vbody, & \omegabodybar &= t_{c} \omegabody,& \fbodybar\extsup &= \frac{1}{N} \fbody\extsup, & \taubodybar\extsup &= \frac{1}{\ell \, N} \taubody\extsup
\end{align*}
\begin{align*}
	 \pbodybar &= \frac{t_{c}}{k \, \ell} \pbody, & \Lbodybar &= \frac{t_{c}}{k \, \ell^2} \Lbody,  &\dragbodybar_{11} &= \frac{1}{\ell \, \eta} \dragbody_{11}, & \dragbodybar_{12} &= \frac{1}{\ell^2 \eta} \dragbody_{12}, & \dragbodybar_{22} &= \frac{1}{\ell^3 \eta} \dragbody_{22}.\end{align*}
The equations \eqref{balance body generic} then become
\begin{equation}
	\label{balance body generic nondim}
	\begin{aligned}
		&\frac{k \, \ell}{t_{c}^2} \left( \dot{\pbodybar} + \omegabodybar \times \pbodybar \right) &= & - \frac{\eta \, \ell^2}{t_{c}} \left( \dragbodybar_{11} \vbodybar + \dragbodybar_{12} \omegabodybar \right) && + N \, \fbodybar\extsup \\
		&\frac{k \, \ell^2}{t_{c}^2} \left( \dot{\Lbodybar} + \omegabodybar \times \Lbodybar \right) &= & - \frac{\eta \,\ell^3}{t_{c}} \left( \dragbodybar_{12}^T \vbodybar + \dragbodybar_{22} \omegabodybar \right) && + \ell \, N \, \taubodybar\extsup.
	\end{aligned}
\end{equation}

We made the assumption that the Reynolds number $\Rey = \rho_{f} \ell^2 / \eta \, t_{c} \ll 1$; accordingly, $\varepsilon := V \Rey / \ell^3$ is small and we can rewrite \eqref{balance body generic nondim} as
\begin{equation}
	\label{balance body nondim eps generic}
	\begin{aligned}
		&\frac{\varepsilon}{1 + \eps_b} \left( \dot{\pbodybar} + \omegabodybar \times \pbodybar \right) &=& - \left( \dragbodybar_{11} \vbodybar + \dragbodybar_{12} \omegabodybar \right) && + \fbodybar\extsup \\
		&\frac{\varepsilon}{1 + \eps_b} \left( \dot{\Lbodybar} + \omegabodybar \times \Lbodybar \right) &=& - \left( \dragbodybar_{12}^T \vbodybar + \dragbodybar_{22} \omegabodybar \right) && + \taubodybar\extsup,
	\end{aligned}
\end{equation}
where $\eps_b = \rho_{f} V/k -1$. Note that the limit $\varepsilon \to 0$ corresponds to the Stokes flow limit $\Rey \to 0$. The limit $\eps \to 0$ could also be satisfied without $\Rey \to 0$, in a thin rod for example. However, hydrodynamic loads defined by equation\re{hydrodynamic loads} require $\Rey \ll 1$ to be valid.

The parameter $\eps_{b}$ is defined such that if the external loads arise solely from buoyancy, the body sinks if $\eps_b<0$, floats if $\eps_b>0$, and is neutrally buoyant if $\eps_b=0$. In principle, the limit $\eps_{b} \to -1$ could be studied, and the inertial effects would then need to be taken into account. In practice however, even for swimmers made of a metallic alloy, i.e. swimmers that have a high density compared to the density of the fluid, the order of magnitude of the swimmer's density $\rho_{s}$ will not exceed about $10 \, \rho_{f}$, so that $1 + \eps_{b}$ is of order $10^{-1}$. Then, one must require $\eps \ll 10^{-1}$ to ensure that $\eps/(1+\eps_{b}) \ll 1$ so that inertial effects can be neglected.

Note that $\frac{\eps}{1+\eps_b} = \frac{k \, N}{l^3 \, \eta^2}$ consistently with the small parameter expressed in~\cite{Gonzalez2004}. In writing it this way, we implicitly aim at comparing the impact of loads arising from gravity and buoyancy in comparison to other external loads.

\subsection{Leading-order Solution}
\label{section: leading order}

The system of differential equation\re{balance body nondim eps generic} is singular. Standard singular perturbation techniques~\citep[]{Hinch1991} are used to find a uniform leading order solution:
\begin{equation}
	\label{leading order sol generic}
	\begin{bmatrix}
		\pbodybar \left( \tbar \right) \\
		\Lbodybar \left( \tbar \right)
	\end{bmatrix}
	= \exp \left( - \left( 1 + \eps_b \right) \frac{\tbar}{\varepsilon} G \right)
	\left(
	\begin{bmatrix}
		\pbodybar_{0} \\
		\Lbodybar_{0}
	\end{bmatrix}
	-
	G^{-1}
	\begin{bmatrix}
		\fbodybar^{(ext)} \left( 0 \right) \\
		\taubodybar^{(ext)} \left( 0 \right)
	\end{bmatrix}
	\right)
	+
	G^{-1}
	\begin{bmatrix}
		\fbodybar^{(ext)} \left( \tbar \right) \\
		\taubodybar^{(ext)} \left( \tbar \right)
	\end{bmatrix},
\end{equation}
where $G := \overline{\Drag} \begin{bsmallmatrix} \Id & 0 \\ 0 & \inertbodybar \end{bsmallmatrix}$, and $\pbodybar_{0}$, $\Lbodybar_{0}$ are the initial conditions to \eqref{balance body nondim eps generic}. This solution is valid as long as the loadings $\fbodybar^{(ext)}$ and $\taubodybar^{(ext)}$ vary slowly compared to the timescale $\varepsilon$ of the initial layer. That is we assume $||d\fbodybar/d\tbar|| \ll \frac{1}{\varepsilon}$ and $||d\taubodybar/\tbar|| \ll \frac{1}{\varepsilon}$.

\subsection{Loads Due to a Uniform Magnetic Field}
\label{section: magnetic loads}

System~\eqref{balance body generic} was studied by \cite{Gonzalez2004} for the particular case where the external loads arise from a constant force, namely the resultant of gravity and buoyancy. Here we treat a complementary case where the external force is $\0$ and the external torque is periodic and has a simple expression: the loads are due to a spatially uniform magnetic field $\Blab$ rotating at constant angular velocity $\alpha$ about a fixed axis of rotation in the lab frame, and the body subjected to it is a permanent magnet so that its magnetic moment $\mlab$ is constant in the body frame. The force and torque are then
\begin{align}
	\label{magnetic loads}
	\flab\msup &= \0, & \taulab\msup &= \mlab \times \Blab.
\end{align}
Without loss of generality, we can choose the lab frame so that the basis vector $\elab_{3}$ corresponds with the axis of rotation of the magnetic field and that $\Blab$ lies in the $\left( \elab_{1}, \elab_{3} \right)$-plane at time $t = 0$, so that $\Blab$ can be explicitly written as
\begin{equation}
	\label{mag field def}
	\Blab = R_{3} \left( \alpha t \right) \Blab_{0}, \qquad \text{where } R_{3} \left( s \right) =
	\begin{bsmallmatrix}
		\cos s & - \sin s & 0 \\
		\sin s & \cos s & 0 \\
		0 & 0 & 1
	\end{bsmallmatrix}
	\text{ and } \Blab_{0} = B \begin{bsmallmatrix} \sin \psi \\ 0 \\ \cos \psi \end{bsmallmatrix},
\end{equation}
where $B$ is the magnitude of the magnetic field, and $\psi$ is the angle between $\Blab$ and its axis of rotation $\elab_{3}$. Note that the lab frame is chosen so that $\elab_{3}$ is the axis of rotation of the magnetic field, and not the direction of gravity.

We non-dimensionalise by substituting $N = m \, B/\ell$, where $m$ is the magnitude of the magnetic moment. The torque due to magnetism is then rewriten in the body frame as
\begin{equation}
	\taubodybar\msup \left( \tbar \right) = \mbodybar \times \Bbodybar \left( \tbar \right) = \mbodybar \times \left( R^T \left( \tbar \right) R_{3} \left( \Ma \, \tbar \right) \begin{bsmallmatrix} \sin \psi \\ 0 \\ \cos \psi \end{bsmallmatrix}\right),
	\label{eqtaubar}
\end{equation}
where all the dependencies in $\tbar$ are written explicitly, and $\Ma$ is the Mason number~\citep{Man2013} given by
\begin{equation}
	\label{Mason number}
	\Ma = \alpha \, t_{c} = \frac{\alpha \, \eta \, \ell^3}{m \, B} \, .
\end{equation}

We define the magnetic frame so that it is aligned with the lab frame at $\tbar =0$ and rotates with the magnetic field. Namely, the magnetic frame is given by the columns of the matrix  $R_{3} \left( \Ma \, t \right)$. We observe in\re{eqtaubar} that the torque $\taubodybar$ applied by the field on the swimmer depends on the  rotation matrix 
\begin{equation}\label{def Q}
Q\left( t \right) := R^T( t ) ~R_{3} ( \Ma \, t )
\end{equation}
 that applies the body frame onto the magnetic frame: $\taubodybar\msup = \mbodybar \times Q \begin{bsmallmatrix} \sin \psi \\ 0 \\ \cos \psi \end{bsmallmatrix}$.

This rotation matrix $Q$ will prove to be more practical than the original matrix $R$. Note that $R$ can be recovered at any point by inverting\re{def Q}. This, in turn, can be substituted in\re{kinematics} to find an equivalent form for Q:
\begin{equation}
	\label{ode Q}
	\dot{Q} = \cross[\left( \Ma \, \ebody_{3} - \omegabody \right)] Q, \qquad Q \left( 0 \right) = Q_{0},
\end{equation}
where $\ebody_{3} = Q \begin{bsmallmatrix} 0 \\ 0 \\ 1 \end{bsmallmatrix}$ is the third basis vector of the lab frame expressed in the body frame.

\subsection{Modelling Gravity and Buoyancy}
\label{section: gravity}

Our study will focus on neutrally buoyant bodies. Nevertheless, this assumption is never strictly achieved in experiments. Therefore, we pause to discuss under which conditions this assumption can be expected to hold.  Let  $-\gbody$ be the force due to gravity (constant in the lab frame) and $\Deltabody$ the misalignment between the centre of mass and the centre of volume (constant in the body frame). We can then write the loads due to gravity and buoyancy as
\begin{align}
	\label{buoyancy loads}
	\fbody\bsup &= \eps_b \gbody, & \taubody\bsup &= \left( 1 + \eps_b \right) \Deltabody \times \gbody.
\end{align}
The corresponding non-dimensionalised quantities in the body frame are $\gbodybar = \frac{1}{k \, g} \gbody$, and $\Deltabodybar = 1/\ell \, \Deltabody$, where $g$ is the gravitational acceleration on earth. Defining the constant $\gamma := \rho_{f} V \, \ell \, g / m \, B$, we rewrite \eqref{balance body nondim eps generic} as
\begin{equation}
	\label{balance body nondim eps}
	\begin{aligned}
		&\frac{\varepsilon}{1 + \eps_b} \left( \dot{\pbodybar} + \omegabodybar \times \pbodybar \right) &=& - \left( \dragbodybar_{11} \vbodybar + \dragbodybar_{12} \omegabodybar \right) &&&& + \gamma \, \frac{\eps_b}{1 + \eps_b} \gbodybar \\
		&\frac{\varepsilon}{1 + \eps_b} \left( \dot{\Lbodybar} + \omegabodybar \times \Lbodybar \right) &=& - \left( \dragbodybar_{12}^T \vbodybar + \dragbodybar_{22} \omegabodybar \right) && + \mbodybar \times \Bbodybar && + \gamma \, \Deltabodybar \times \gbodybar.
	\end{aligned}
\end{equation}

In the following, we focus on the case of a neutrally buoyant body whose centre of mass and centre of buoyancy coincide, that is $\eps_b = 0$ and $\Deltabodybar = \0$. These strict equalities cannot be achieved experimentally; however we expect our conclusions to be a good approximation when $|\eps_b|, \left| \Deltabodybar \right| \ll 1/\gamma$. See~\cite{Ruegg-Reymond2019} for a study of cases where these hypotheses are relaxed.

\subsection{Governing Equations}
\label{section: get equations}

Assuming that the swimmer is neutrally buoyant, we substitute $\eps_{b} = 0$ and $\Deltabodybar = \0$ in \eqref{balance body nondim eps}. Gathering~(\ref{kinematics}, \ref{velocity momentum}, \ref{ode Q}) and dropping the overbars, we obtain the full system of differential equations for our problem
\begin{eqnarray}
	\label{balance body nondim no buoyancy}
	&	\begin{aligned}
 \dot{\xlab} &= \vlab,\\
 \dot Q &= [\Ma\, \ebody_3 -\omegabody]^\times \,Q,\\
		\varepsilon \left( \dot{\pbody} + \omegabody \times \pbody \right) &= - \left( \dragbody_{11} \vbody + \dragbody_{12} \omegabody \right) \\
		\varepsilon \, \left( \dot{\Lbody} + \omegabody \times \Lbody \right) &= - \left( \dragbody_{12}^T \vbody + \dragbody_{22} \omegabody \right) + \mbody \times \Bbody,
	\end{aligned}
\end{eqnarray}
together with the relations 
\begin{equation} \label{balance body nondim no buoyancy constitutive law}
\pbody = \vbody,\qquad\qquad	\Lbody = \inertbody \omegabody,\qquad\qquad\Bbody \left( t \right) = Q \left( t \right) \begin{bsmallmatrix}
		\sin \psi \\ 0 \\ \cos \psi
	\end{bsmallmatrix}.
\end{equation}

\subsection{Long-time Behaviour}
\label{section: long-time behav}

In the particular case of the problem\ret{balance body nondim no buoyancy}{balance body nondim no buoyancy constitutive law}, the solution \eqref{leading order sol generic} becomes \begin{equation}
	\label{leading order sol}
	\begin{bmatrix}
		\pbody \left( t \right) \\
		\Lbody \left( t \right)
	\end{bmatrix}
	= \exp \left( -\frac{t}{\varepsilon} G \right)
	\begin{bmatrix}
		\pbody_{0} -  \mobbody_{12} \cross[\mbody] \Bbody \left( 0 \right) \\
		\Lbody_{0} - \inertbody \mobbody_{22} \cross[\mbody] \Bbody \left( 0 \right)
	\end{bmatrix}
	+
	\begin{bmatrix}
		\mobbody_{12} \cross[\mbody] \Bbody \left( t \right) \\
		\inertbody \mobbody_{22} \cross[\mbody] \Bbody \left( t \right)
	\end{bmatrix}
\end{equation}
where $\mobbody_{12}$ and $\mobbody_{22}$ are 3-by-3 blocks of the mobility matrix
\begin{align*}
	\Mob = \Drag^{-1} = 
	\begin{bmatrix}
		\mobbody_{11} & \mobbody_{12} \\
		\mobbody_{12}^T & \mobbody_{22}
	\end{bmatrix} \, ,
\end{align*}
which is constant in the body frame.

In order for \eqref{leading order sol} to be a valid leading order solution, we need $\Ma \ll 1/\varepsilon$. We expect that it represents accurately the dynamics of our system in the outer layer, i.e. for $t \gg \varepsilon$, in the limit $\varepsilon \to 0$. In the inner layer however, the leading-order solution may be inaccurate for any $\varepsilon\neq0$ (see \cite{Gonzalez2004}  and \cite{Kim2013}). 

In the outer layer, the first term in\re{leading order sol} can be neglected and substituting the first two equations of \re{balance body nondim no buoyancy constitutive law} therein, the problem\ret{balance body nondim no buoyancy}{balance body nondim no buoyancy constitutive law} becomes
\begin{equation}
	\label{outer layer system}
	\begin{aligned}
		\dot{\xlab} &= \vlab,\\
 \dot Q &= [\Ma\, \ebody_3 -\omegabody]^\times \,Q,\\
 \end{aligned} 
 \end{equation}
 where 
\begin{equation}
	\label{cl}
	\begin{aligned}
		& \vbody \left( t \right) = \mobbody_{12} \cross[\mbody] \Bbody \left( t \right), &&\qquad&
		& \omegabody \left( t \right) = \mobbody_{22} \cross[\mbody] \Bbody \left( t \right), \\
		& \Bbody \left( t \right) = Q\left( t \right) \begin{bsmallmatrix} \sin \psi \\ 0 \\ \cos \psi \end{bsmallmatrix}, &&\qquad&
		& \ebody_{3} \left( t \right) = Q\left( t \right) \begin{bsmallmatrix} 0 \\ 0 \\ 1 \end{bsmallmatrix}.
	\end{aligned}
\end{equation}

Finally, substituting\re{cl} in\re{outer layer system}, we find the following closed form system of ODEs
\begin{align}
	&\dot{\xlab} =R_3(\Ma\,t) \,Q(t)\,  \mobbody_{12}\, \cross[\mbody] \, Q\, \begin{bsmallmatrix} \sin \psi \\ 0 \\ \cos \psi \end{bsmallmatrix}, &&\qquad&
	&\xlab(0) = \xlab_0, \label{ode x} \\
	&\dot{Q} = \cross[\left( \Ma \, Q
	\begin{bsmallmatrix}
		0 \\ 0 \\ 1
	\end{bsmallmatrix}
	 - \Pbody \, Q
	 \begin{bsmallmatrix}
		\sin \psi \\ 0 \\ \cos \psi
	\end{bsmallmatrix}
	 \right)] Q ,&&\qquad&
	 &Q \left( 0 \right) = Q_{0},	\label{ode Q closed}
\end{align}
where $\Pbody := \mobbody_{22} \cross[\mbody]$ is a constant matrix that only depends on the shape and the magnetisation of the swimmer. Propositions 3.1 and 3.2 of \cite{Gonzalez2004} can be straightforwardly adapted to show that the solutions of\ret{ode x}{ode Q closed} completely determine the long-time behaviour of the leading-order solution of\ret{balance body nondim no buoyancy}{balance body nondim no buoyancy constitutive law}. 

Much of our subsequent analysis will consist in studying the set of solutions of\ret{ode x}{ode Q closed}. We readily note that\re{ode Q closed} decouples from\re{ode x} so that for any solution of the $\dot{Q}$ equation, $\xlab$ can be recovered by quadrature. Accordingly most of what follows will be  focused on\re{ode Q closed}.

\section{Analysis of Governing Equations}
\label{section: analysis}

In this section, we analyse the dynamical system\ret{ode x}{ode Q closed}. In particular, we study how the number of relative equilibria of the system and their stability depend on two parameters that can be changed during an experiment: the Mason number $\Ma$ given by\re{Mason number} and the angle $\psi$ between the magnetic field and its axis of rotation. The equation\re{ode Q closed} is studied in sections~\ref{sec-releq},~\ref{section: symmetry} and~\ref{section: stability}. The swimmers trajectories in the lab frame are then recovered in section~\ref{sec-traj} by studying\re{ode x}.


\subsection{The set of relative equilibria of a given magnetic swimmer}
\label{sec-releq}

The swimmer is at \emph{relative equilibrium} when $Q$ is a steady state solution of the differential equation in\re{ode Q closed}: the orientation of the swimmer is fixed in the magnetic frame. Accordingly, relative equilibria occur if and only if
\begin{equation} 
	\label{eq cond Q}
	\Ma \, Q \, \begin{bsmallmatrix} 0 \\ 0 \\ 1 \end{bsmallmatrix} = \Pbody \, Q \, \begin{bsmallmatrix} \sin \psi \\ 0 \\ \cos \psi \end{bsmallmatrix}.
\end{equation}
  Instead of working directly with the rotation matrix $Q$, we will work with the body frame components of two vectors: $\ebody_3$ and $\Bbody$. In this picture, finding relative equilibria (i.e. $Q$ respecting\re{eq cond Q}) amounts to finding pairs $\left( \ebody_{3}, \Bbody \right)$ such that
\begin{align}
	\label{equilibrium conditions}
	\Ma \, \ebody_{3} &= \Pbody \, \Bbody, 
	& \ebody_{3} \cdot \Bbody &= \cos \psi,
	 &\ebody_3 \cdot \ebody_3 &=1, 
	 & \Bbody\cdot \Bbody &=1,
\end{align}
for given parameters $\Ma$ and $\psi$. Note that an equivalent set of equations can be found in~\cite{Meshkati2014,Fu2015} where they were solved numerically. The version that we propose here enables the determination of all solutions analytically.

Although experimentally, one wants to fix $\Ma$ and $\psi$ and then find all the $(\ebody_3, \Bbody)$ solutions of~\eqref{equilibrium conditions} for these parameters, it is mathematically expedient to consider~\eqref{equilibrium conditions} as a system of 6 equations of the 8 unknown $(\Ma,\cos\psi,\,\ebody_3,\,\Bbody)$. We will show that the complete set of solutions is then a two dimensional surface that can be parameterised by a single map. That is, each pair of values of the parameters to be described here under corresponds to a single equilibrium of the system. This is in contrast with the experimental viewpoint where for a given $\Ma$ and $\psi$ there can be up to 8 different solutions of~\eqref{equilibrium conditions}.
%


We first consider the singular value decomposition of the matrix $\Pbody = \mobbody_{22} \cross[\mbody]$. Let $\sigma_{0} = 0< \sigma_2<\sigma_{1}$ be the singular values of $\Pbody$ with corresponding right-singular vectors $\betabody_{0} = \mbody$, $\betabody_{1}$, and $\betabody_{2}$ and left-singular vectors $\etabody_{0} = \mobbody_{22}^{-1} \mbody / \left\| \mobbody_{22}^{-1} \mbody \right\|$, $\etabody_{1}$, and $\etabody_{2}$. Recall that we have the relations
\begin{align}
	\Pbody \betabody_{i} &= \sigma_{i} \etabody_{i}, & \Pbody^T \etabody_{i} = \sigma_{i} \betabody_{i},\label{singvecttheory}
\end{align}
and that the two sets of singular vectors can be chosen so that they each form a right-handed and orthonormal basis. 

The first equilibrium condition in \eqref{equilibrium conditions} constrains $\ebody_{3}$ to lie in the span of $\Pbody$, that is in the $\left( \etabody_{1}, \etabody_{2} \right)$-plane. Therefore, an angle $\theta \in \left( -\pi, \pi \right]$ exists such that
\begin{equation}
	\label{e3 param}
	\ebody_{3} = \cos \theta \, \etabody_{1} + \sin \theta \, \etabody_{2}.
\end{equation}

Next, we expand $\Bbody$ in the $\left\{ \betabody_{i} \right\}$-basis as $\Bbody = \cos \phi \, \betabody_{0} + \sin \phi \left( \cos \xi \, \betabody_{1} + \sin \xi \, \betabody_{2} \right)$ for $\phi \in \left[ 0 , \pi \right]$ and $\xi \in \left[ 0, 2 \pi \right)$. Substituting in\re{equilibrium conditions} and using\re{singvecttheory}, we find that
\begin{itemize}
\item ~there is a one-to-one correspondence between $\xi$ and $\theta$ so that $\Bbody$ rewrites as
\begin{equation}
	\label{B param}
	\Bbody = \cos \phi \, \betabody_{0} + \sin \phi \left( \frac{\cos^2 \theta}{\sigma_{1}^2} + \frac{\sin^2 \theta}{\sigma_{2}^2} \right)^{-1/2} \, \left( \frac{\cos \theta}{\sigma_{1}} \, \betabody_{1} + \frac{\sin \theta}{\sigma_{2}} \, \betabody_{2} \right),
\end{equation}
and
\item ~$\Ma$ and $\cos \psi$ are also fixed by $\theta$ and $\phi$:
\end{itemize}
\begin{subequations}
	\label{param a psi}
	\begin{align}
		\Ma =& \sin \phi \left( \frac{\cos^2 \theta}{\sigma_{1}^2} + \frac{\sin^2 \theta}{\sigma_{2}^2} \right)^{-1/2}, \label{param a} \\
		\cos \psi =& \cos \phi \left( c_{01} \, \frac{\cos \theta}{\sigma_{1}} + c_{02} \, \frac{\sin \theta}{\sigma_{2}} \right) \label{param psi} \\
		&+ \sin \phi \left( \frac{\cos^2 \theta}{\sigma_{1}^2} + \frac{\sin^2 \theta}{\sigma_{2}^2} \right)^{-1/2} \, \left( c_{11} \, \left( \frac{\cos^2 \theta}{\sigma_{1}^2} - \frac{\sin^2 \theta}{\sigma_{2}^2} \right) + c_{12} \, \frac{\cos \theta \, \sin \theta}{\sigma_{1} \, \sigma_{2}} \right) \, , \nonumber
	\end{align}
\end{subequations}
where $c_{01} = \betabody_{0} \cdot \Pbody \, \betabody_{1}$,  $c_{02} = \betabody_{0} \cdot \Pbody \ \betabody_{2}$, $c_{11} = \betabody_{1} \cdot \Pbody \, \betabody_{1} = - \betabody_{2} \cdot \Pbody \, \betabody_{2}$, and $c_{12} = \betabody_{1} \cdot (\Pbody + \Pbody^T) \, \betabody_{2}$.

As a result, the angles $\theta$ and $\phi$ are smooth coordinates on the set of relative equilibria. That is a value $\left( \theta, \phi \right)$ is in one-to-one correspondence with a unique relative equilibrium.

\begin{figure}
	\includegraphics{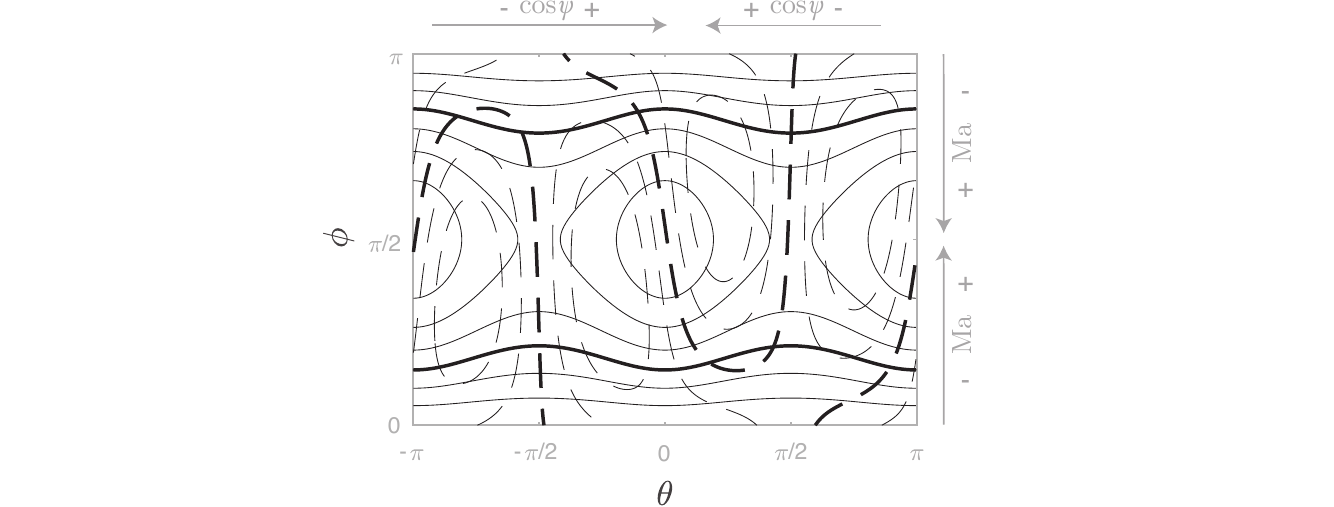}
	\caption[]{Level sets of the Mason number (solid lines) and $\cos \psi$ (dashed lines) as functions of the two mapping parameters $\theta$ and $\phi$.  Note that the figure can be smoothly wrapped on a cylinder by identifying the borders at $\theta = \pm \pi$. The data shown here correspond to swimmer A (cf. section~\ref{section: numerics}).
	}
	\label{fig: contours Ma cospsi}
\end{figure}

The level sets of $\Ma$ and $\cos \psi$ are shown in figure~\ref{fig: contours Ma cospsi} for the helical swimmer~A described in section~\ref{section: numerics} and shown in figure~\ref{fig: swimmers}. Figure~\ref{fig: contours Ma cospsi} can be used to find all the equilibria occurring for any given values of the experimental parameters. The arrows indicating increasing values of $\cos \psi$ are valid on the upper border of the figure. The bold curves highlight the particular level sets corresponding to $\Ma=0.015$ and $\cos\psi= 0.01$. For these parameter values, the full and dashed lines have 8 intersections corresponding to 8 different relative equilibria at these values of the parameters. Once the pair ($\theta,\phi$) is found for a particular equilibrium, the orientation of the body can be inferred by computing $\ebody_3$ and $\Bbody$ from~\eqref{e3 param} and \eqref{B param}. Note that this map also proves that there are relative equilibria for all orientations of $\Bbody$ (not all of them stable, more on this in section~\ref{section: stability}) but not for all orientations of $\ebody_3$: at equilibrium, the rotation axis is in a plane orthogonal to $\mobbody_{22}^{-1} \mbody$.

The maximal $\Ma$ for which a given swimmer has any relative equilibria is equal to the largest singular value $\sigma_{1}$ of the $\Pbody$ matrix (see equation~\eqref{param a}). It is straightforward to show that for a swimmer of a given shape,  the maximal value of $\sigma_{1}$ is obtained by magnetising the body in a direction perpendicular to the "easy rotation axis" -- i.e. perpendicular to the eigenvector of $\mobbody_{22}$ corresponding to its largest eigenvalue. 

Finally, note that by definition~\eqref{e3 param}, all functions of $\theta$ can be extended to the real line by imposing periodicity modulo $2 \pi$. In particular, it will help to wrap figure~\ref{fig: contours Ma cospsi} on a cylinder by identifying the borders at $\theta = \pm \pi$.

\subsection{Symmetry of the rotational dynamic}
\label{section: symmetry}

The solutions of the ODE~\eqref{ode Q closed} come in symmetric pairs. Indeed, if the function $Q \left( t \right)$ is a solution, then so is 
$$\breve{Q} \left( t \right) := Q \left( -t \right) \begin{bsmallmatrix} -1 & \phantom{-}0 & \phantom{-}0 \\ \phantom{-}0 & \phantom{-}1 & \phantom{-}0 \\ \phantom{-}0 & \phantom{-}0 & -1 \end{bsmallmatrix} =: Q \left( -t \right) R_{2} \left( \pi \right).$$

With obvious notations, one finds that in the body frame  $\breve{\Bbody} \left( t \right) = - \Bbody \left( -t \right)$, $\breve{\ebody}_{3} \left( t \right) = - \ebody_{3} \left( -t \right)$, $\breve{\vbody} \left( t \right) = - \vbody \left( -t \right)$ and $\breve{\omegabody} \left( t \right) = - \omegabody \left( -t \right)$.

In the magnetic frame, the equations of motion are autonomous. At any instant, if the body is oriented according to $Q^T$ (w.r.t. the magnetic frame) and has angular velocity $\ulab$ with respect to the magnetic frame, then the symmetric solution is rotated by $\pi$ through the axis perpendicular to both $\elab_3$ and $\Blab$ and the angular velocity $\breve{ \ulab} $ is the reflection of $\ulab$ through the plane containing both $\elab_3$ and $\Blab$.

At relative equilibria, the body is locked in the magnetic frame ($\ulab= \boldsymbol 0$) and symmetric solutions are obtained by rotations through $\pi$ with respect to the axis perpendicular to both $\elab_3$ and $\Blab$. Accordingly, one finds 
\begin{align} 
\label{symeq}
\breve \phi &= \pi - \phi.&
\breve \theta &=(2\pi+\theta \mod 2\pi)-\pi.
\end{align} 


The corresponding symmetry of figure~\ref{fig: contours Ma cospsi} consists in a reflection through the line at $\phi = \pi/2$ followed by a horizontal translation by $\pi$ (making use of the $\theta$-periodicity).

\subsection{Stability}
\label{section: stability}

\begin{figure}[b!]
	\centering
	\includegraphics{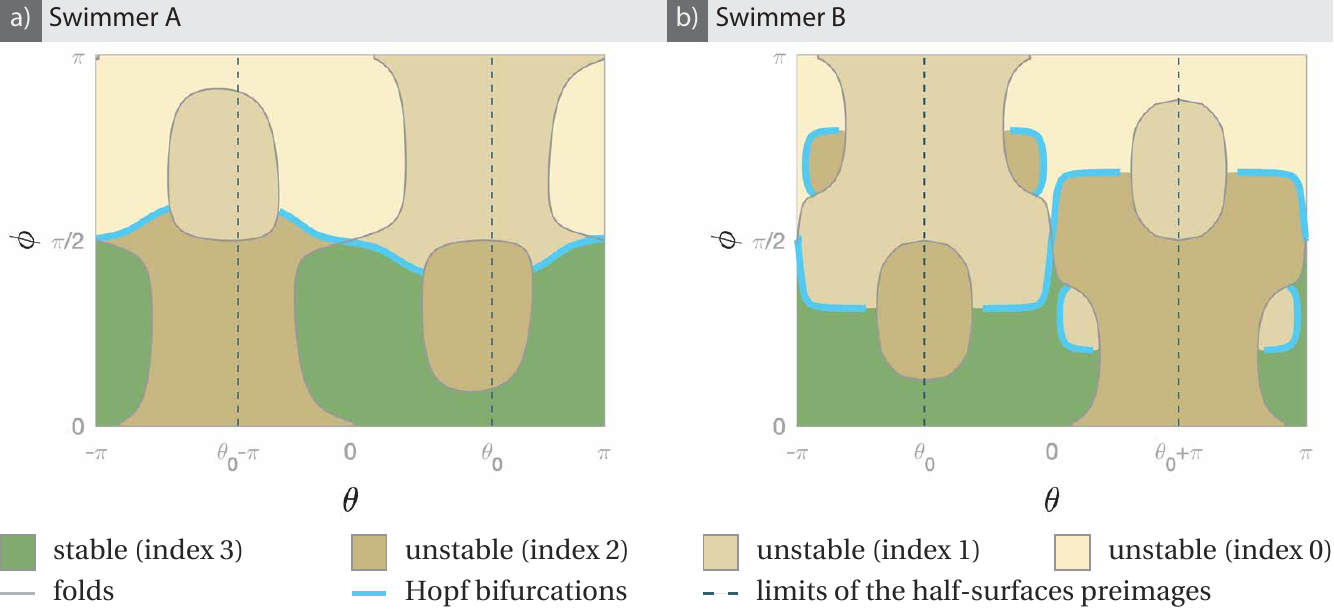}
	{\phantomsubcaption\label{subfig: stabchart a}}
	{\phantomsubcaption\label{subfig: stabchart b}}
	\caption{Stability index as a function of the two mapping parameters $\theta$ and $\phi$, with folds and Hopfs bifurcations at the boundaries between regions with different stability indices. The data shown here corresponds to the  swimmers A  and B described in section~\ref{section: numerics}.}
	\label{fig: stabchart}
\end{figure}

The linear stability of relative equilibria of \eqref{outer layer system} is given by the sign of the real part of the eigenvalues of the linearised dynamics matrix
\begin{equation}
A = \Pbody \cross[\Bbody] - \Ma \cross[\ebody_{3}], \label{StabMat}
\end{equation}
where $\Bbody$ and $\ebody_3$ are taken at that equilibrium. They are therefore explicit functions of $\theta$ and $\phi$. Each equilibrium is thus associated with an index defined as the number of eigenvalues of $A$ with strictly negative real parts (counting multiplicities). An equilibrium is therefore unstable whenever its index is strictly less than $3$. If the index is 3, it is stable (see figure~\ref{fig: stabchart}). 

 We remark that the linearised dynamics matrices $A$ of symmetric pairs of relative equilibria (in the sense of section~\ref{section: symmetry}) have opposite eigenvalues. So if one equilibrium of the pair is stable, then the other is necessarily unstable with index 0. 

As we move on the set of relative equilibria, stability changes only at bifurcation points of the dynamical system \eqref{ode Q closed}~\citep[]{Kuznetsov2004}. In this system almost all  bifurcations are either folds or Hopf bifurcations. 

\subsubsection{Folds}

Fold points are generically points where one real eigenvalue of $A$ is $0$, so $\det A = 0$. At a fold, the linearised system admits a line of equilibria. This occurs when the parameterisation of the $(\Ma,\cos\psi)$ plane by $(\theta,\phi)$ looses regularity, that is when the Jacobian
\begin{align} \label{JacobianFold}
\left | 
\frac {\partial\, ( \Ma,\cos \psi)}{\partial \,( \theta,\phi)}\right | = 0,
\end{align}
vanishes.
This condition makes them easy to track since we have the explicit parameterisation~\eqref{param a psi}. 

Finally,  folds also correspond to the values of the experimental parameters $\Ma$ and $\cos \psi$ for which the number of equilibria changes. Therefore folds occur when the level sets of $\Ma$ and $\cos \psi$ in figure~\ref{fig: contours Ma cospsi} are tangent\footnote{For fixed values of $\Ma$ and $\cos\psi$, points of tangency of the level curves go by pair because of the symmetry described in section~\ref{section: symmetry}.}.

%


\subsubsection{Hopf Bifurcations}
\label{section: hopf bifurcations}

Hopf bifurcations are points where $A$ has a pair of purely imaginary conjugate eigenvalues $\pm i \, \lambda_{I}$. Accordingly, relative equilibria on both sides of the Hopf bifurcation have indices differing by $2$. 

To find the Hopf bifurcations of~\eqref{ode Q closed}, we apply the technique described in~\cite{Kuznetsov2004} and compute the bialternate product
\begin{align} 
2 \, A \odot I =
\begin{pmatrix} 
a_{11} + a_{22}  & a_{23} & - a_{13}  \\
a_{32}  & a_{33}+ a_{11}   & a_{12}\\
-a_{31}  & a_{21} & a_{22}  + a_{33}
\end{pmatrix} ,
\end{align}  where $I$ is the identity matrix in $\mathbb R^{3\times 3}$ and $a_{ij}$ are the elements of $A$. By construction, the eigenvalues of $2 A\odot I$ are the three sums of pairs of eigenvalues of $A$. The Hopf bifurcations are therefore obtained as the zero-level curve of $\det (2 \, A \odot I)$. Note that $\det (2 \, A \odot I)$ also cancels when $A$ has two opposite real eigenvalues; these are removed in post-processing.

Hopf bifurcations also give rise to branches of periodic orbits. We used numerical continuation techniques to explore them and the results are exposed in section~\ref{section: periodic}.

\subsection{Translational trajectories corresponding to specific relative equilibria}
\label{sec-traj}

By definition, the trajectory corresponding to a relative equilibrium has constant linear and angular velocities $\vbody$ and $\omegabody$ in the body frame and is therefore helical. The axis of the helix is given by $\omegabody$ while its pitch and radius are given by~\citep{Crenshaw1993}.
\begin{equation}
		p = \frac{2 \pi}{\Ma^2} \, \omegabody \cdot \vbody \quad \text{and} \quad r = \frac{1}{\Ma^2} \left| \omegabody \times \vbody \right|. \label{heltheo}
\end{equation}
Substituting the equilibrium condition\re{equilibrium conditions} in the equation for the outer layer behaviour\re{cl} yields \begin{align} 
\vbody & = \Ma\, \mobbody_{12} \mobbody_{22}^{-1} \, \ebody_3, & \omegabody &= \Ma\, \ebody_3.\label{vomst}
\end{align} 
Further substituting\re{vomst} in\re{heltheo}, we obtain
\begin{equation*}
		p = 2 \pi \, \ebody_{3} \cdot \mobbody_{12} \mobbody_{22}^{-1} \ebody_{3} \quad \text{and} \quad r = \left| \ebody_{3} \times \mobbody_{12} \mobbody_{22}^{-1} \ebody_{3} \right|.
\end{equation*}
The swimmer evolves along a helix the axis of which is parallel to the axis of rotation of the external field (see eq.~\eqref{vomst}). The so called chirality matrix $\textit{Ch} =\mobbody_{12}\mobbody_{22}^{-1}$ is studied in details in~\cite{Morozov2017} including a discussion of how the off-diagonal terms of $\textit{Ch}$ allow some non-chiral objects to nevertheless transform their rotational dynamic into translational motion.



\subsubsection{Optimal $\Ma$ and $\psi$ for maximal axial velocity of a magnetic swimmer}

\begin{figure}[b!]
	\includegraphics{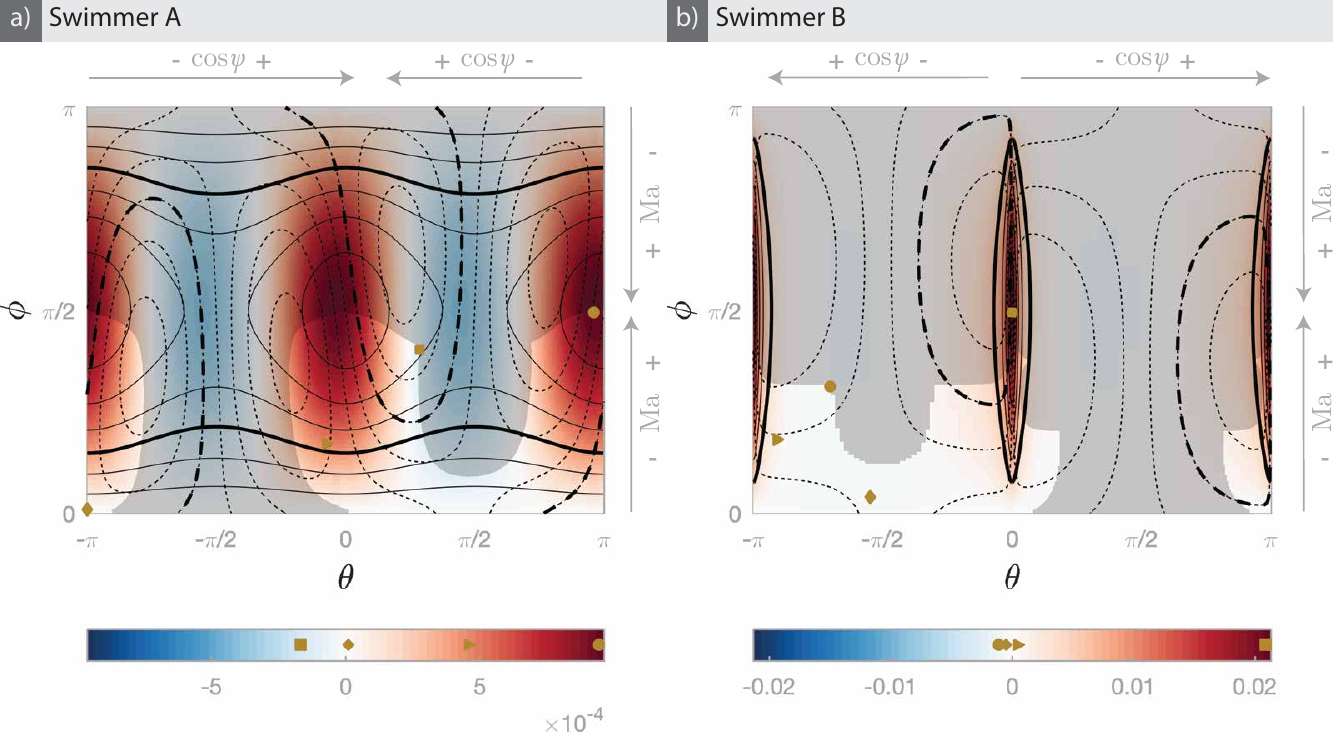}
	{\phantomsubcaption\label{subfig: axvel a}}
	{\phantomsubcaption\label{subfig: axvel b}}
	\caption{Axial velocity in function of the two parameters $\theta$ and $\phi$. The level sets of the Mason number (solid lines) and $\cos \psi$ (dashed lines) are also shown (cf.~fig.~\ref{fig: contours Ma cospsi}). The shaded area comprises only unstable relative equilibria. The data shown here correspond to swimmers~A and B (cf. section~\ref{section: numerics}). In panel~\subref{subfig: axvel b}, bold curves correspond to $\Ma = 0.2198$ and $\cos\psi = -0.3989$.
	}
	\label{fig: contours Ma cospsi axvel}
\end{figure}

The axial velocity $\vax$ is the component of the velocity along the axis $\ebody_3$ of rotation of the magnetic field: \begin{equation}
	\label{ax vel}
	\vax = \Ma \, \ebody_{3} \cdot \mobbody_{12} \, \mobbody_{22}^{-1} \, \ebody_{3}.
\end{equation}

Substituting the parametrisation\ret{e3 param}{param a psi} to express $\Ma$ and $\ebody_{3}$ in terms of $\theta$ and $\phi$ allows to find the maximal axial velocity admissible for relative equilibrium trajectories as
\begin{equation}
	\vax(\theta,\phi) =  \sin \phi \left( \frac{\cos^2 \theta}{\sigma_{1}^2} + \frac{\sin^2 \theta}{\sigma_{2}^2} \right)^{-1/2} \, \left( \cos \theta \, \etabody_{1} + \sin \theta \, \etabody_{2} \right) \cdot \mobbody_{12} \, \mobbody_{22}^{-1} \, (\cos \theta \, \etabody_{1} + \sin \theta \, \etabody_{2}) \, .\label{vax}
\end{equation}

Figure~\ref{fig: contours Ma cospsi axvel} shows the function $\vax$~\eqref{vax} as a colour plot. The greyed area indicates unstable equilibria. For optimal handling, the function $\vax$ may be maximised (in absolute value) with  respect to $\theta$ and $\phi$. The optimal handling parameters $\Ma$ and $\cos \psi$ can then be read from~\eqref{param a psi}.

Note that although it is clear from~\eqref{vax} that the maxima will be reached at $\phi = \pi/2$, which corresponds to maximal $\Ma$ for a given $\theta$ (cf. parametrisation\re{param a}), it is in general false that increasing $\Ma$ increases the axial velocity; see section~\ref{sec: comparison} for an example. 


\subsubsection{Optimal magnetisation for large axial velocity}
\label{section: opt m}
To optimally magnetise a swimmer of a given shape to achieve maximal axial velocity, we make use of $\Ma \, \ebody_{3} = \Pbody \, \Bbody = \mobbody_{22} \, (\mbody \times \Bbody)$ to obtain
\begin{equation*}
	\vax = \frac{(\mbody \times \Bbody) \cdot \mobbody_{22} \, \mobbody_{12} \, (\mbody \times \Bbody)}{|\mobbody_{22} \, (\mbody \times \Bbody)|} \, .
\end{equation*}
Noting that $\mbody \times \Bbody = \sin \phi \, \nbody$ for some unit vector $\nbody$, we consider the expression \begin{equation}
	\label{ax vel n}
\sin \phi \, \frac{\nbody \cdot \mobbody_{22} \, \mobbody_{12} \, \nbody}{|\mobbody_{22} \, \nbody|}.
\end{equation}
Let $\nbody^\star$ be the absolute maximizer of
\begin{equation}
	\label{ax vel to opt}
	f(\nbody) = \frac{|\nbody \cdot \mobbody_{22} \, \mobbody_{12} \, \nbody|}{|\mobbody_{22} \, \nbody|}
\end{equation}
over the unit vectors.

Next, we show that optimal magnetisation is achieved by picking $\mbody$ orthogonal to $\nbody^\star$. Note that while it is clear that~\eqref{ax vel n} is optimised by picking $\phi=\pi/2$ and $\nbody=\nbody^\star$, at this point, $\vax^\star = f(\nbody^\star)$ is an upper bound on the maximal axial velocity of a swimmer of a given shape. We still need to show that there exists a $\theta^\star$ such that $\mbody \times \Bbody(\theta^\star,\pi/2) = \nbody^\star$.  To this intent, remember from the previous section that for a given magnetisation, optimal axial velocity is indeed obtained at $\phi = \pi/2$. In this case, by varying $\theta$, we can pick $\Bbody$ to be any unit vector that is perpendicular to $\mbody$. In particular, it is possible to choose $\theta$ such that $\Bbody (\theta,\pi/2) = \pm \nbody^\star \times \mbody$: the optimal velocity $\vax^\star$ is reached for those values of~$\theta$.

Note that we are left with a one parameter family of magnetisations for which it is possible to find an equilibrium with $\vax =\vax^\star$. This degree of freedom may be used to try to bring the optimal equilibrium close to the stable region ($\phi = \pi/2$ is itself expected to be neutrally stable at best). 
 
Let $x \in \mathbb S^1$ parametrise the unit circle in the plane $\nbody^{\star\,\bot}$ as $\{ \mbody(x) : x \in \mathbb S^1 \}$. For each $x$, compute the eigenvalues of the matrix\re{StabMat}
\begin{equation*}
	A(x) = \mobbody_{22} \, \cross[\mbody(x)] \, \cross[\big(\nbody^\star \times \mbody(x) \big)] - \cross[(\mobbody_{22} \, \nbody^\star)] \, .
\end{equation*}
Look for $x$ such that $A(x)$ has a pair of purely imaginary eigenvalues. If such an $x$ exists, set the swimmer's magnetic moment direction as $\mbody^\star = \mbody(x)$. Then the two equilibria $(\Bbody, \ebody_{3}) = \pm \left( \nbody^\star \times \mbody^\star, \frac{\mobbody_{22} \, \nbody^\star}{|\mobbody_{22} \, \nbody^\star|} \right)$ maximise the axial velocity with parameter values $\Ma = |\mobbody_{22} \, \nbody^\star|$ and $\cos \psi = \big(\nbody^\star \times \mbody^\star \big) \times \frac{\mobbody_{22} \, \nbody^\star}{|\mobbody_{22} \, \nbody^\star|}$, and at least one of them has stable relative equilibria in its neighbourhood.

Indeed the two pairs $(\Bbody, \ebody_{3}) = \pm \left( \nbody^\star \times \mbody(x), \frac{\mobbody_{22} \, \nbody^\star}{|\mobbody_{22} \, \nbody^\star|} \right)$ are then Hopf bifurcations of the system, and their respective stability matrices are opposite. So the unique real eigenvalue of one of the stability matrix is non-positive, while the pair of purely imaginary eigenvalues have zero real part. Therefore, at least one of these two Hopf bifurcation has stable relative equilibria in its neighbourhood.

 An example where $\mbody$ is chosen as to maximise the axial velocity on stable equilibria is presented in section~\ref{section: numerics}. 


\subsection{There are almost always 0, 4, or 8 relative equilibria}
\label{section: surf charac}

\begin{figure}[b!]
	\centering
	\includegraphics{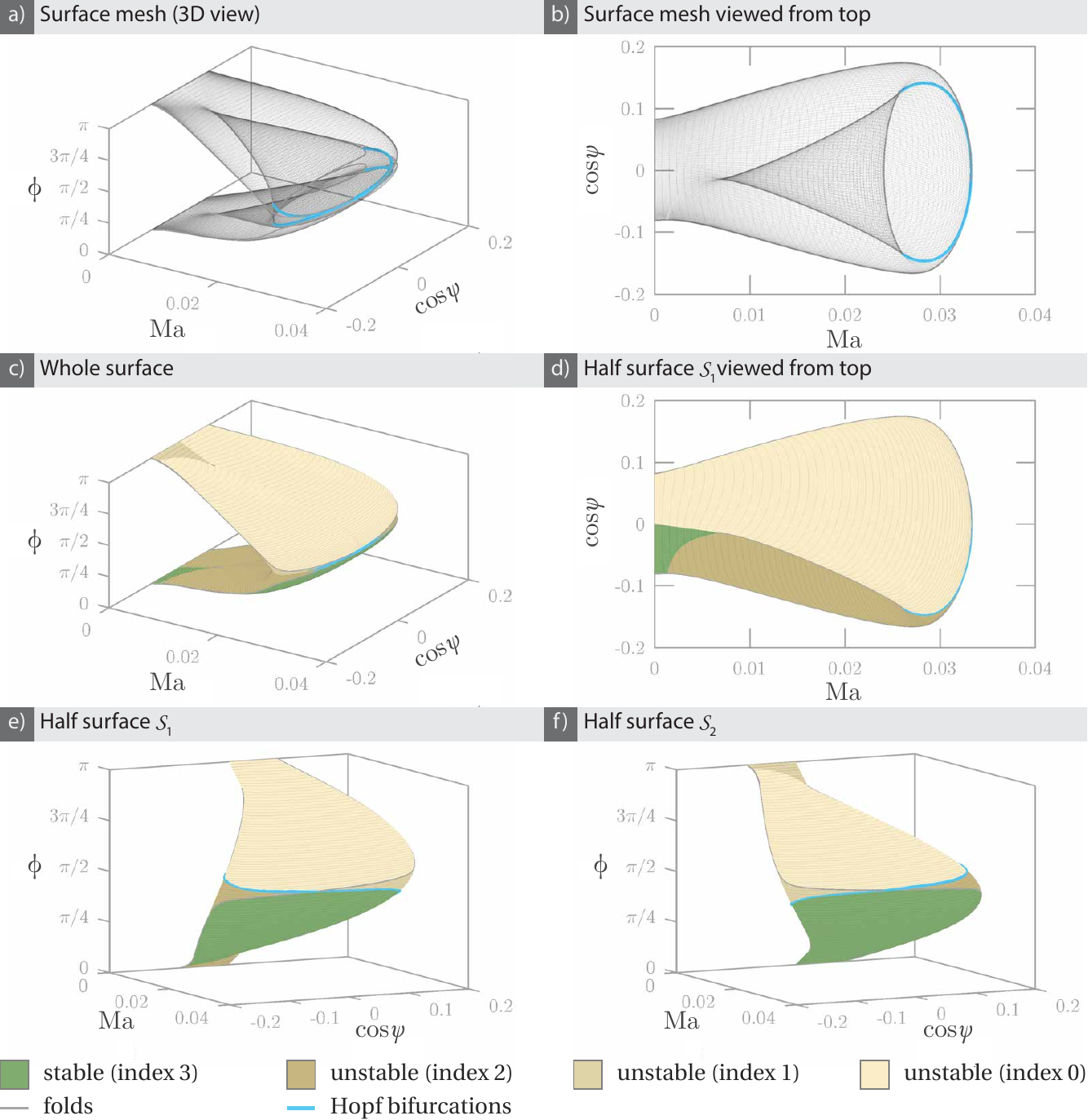}
	{\phantomsubcaption
	\label{subfig: swimmer A mesh}}
	{\phantomsubcaption
	\label{subfig: swimmer A mesh top}}
	{\phantomsubcaption
	\label{subfig: swimmer A surf}}
	{\phantomsubcaption
	\label{subfig: swimmer A surf top}}
	{\phantomsubcaption
	\label{subfig: swimmer A surf 1/2}}
	{\phantomsubcaption
	\label{subfig: swimmer A surf 2/2}}
	\caption[]{Surface $\mathcal S$ parametrised by eq.\re{surface param}, providing a visualisation of the set of relative equilibria of a particular swimmer in the space of coordinates $(\Ma,\cos\psi, \phi)$. The surface $\mathcal S$ is split into two symmetric half surfaces $\mathcal S_{1}$ and $\mathcal S_{2}$ (cf. eq.\re{half surfs}). The data shown here corresponds to swimmer A (cf. section~\ref{section: numerics}).}
	\label{fig: swimmer A eq}
\end{figure}

We have seen is section~\ref{sec-releq} that the set of relative equilibria forms a 2-dimensional surface in the eight dimensional space $(\Ma,\cos\psi,\,\ebody_3,\,\Bbody)$. To visualise it we project it in $\mathbb R^3$ by considering the function
\begin{equation}\label{surface param}
\Sigma :(-\pi,\pi] \times [0,\pi]\mapsto\mathbb R^3: ( \theta, \phi) \mapsto \Big( \Ma \left( \theta, \phi \right), \cos \psi \left( \theta, \phi \right), \phi \Big),
\end{equation}
and defining the surface $\mathcal S$ as the image set of $\Sigma$. Furthermore $\Sigma$ is extended to $\mathbb R\times [0,\pi]$ by asking that it be $2\pi$-periodic in $\theta$.

The graph of $\Sigma$, or equivalently the surface $\mathcal S$, is shown in figure~\ref{fig: swimmer A eq} with the colour encoding the stability index (green is stable). Close inspection of equations~\eqref{param a psi} reveals that $\Sigma$ is not injective (see appendix~\ref{appendix: intersections} for details). One of the self-intersection of $\mathcal S$ is found along the coordinate line $(\theta_0,\phi)$ with $\theta_0 = \arctan \left( -c_{01} \, \sigma_2 / (c_{02} \, \sigma_{1}) \right)$ where $\Sigma \left( \theta_{0}, \phi \right) = \Sigma \left( \theta_{0} + \pi, \phi \right)$ for all $\phi$. Furthermore, the two restrictions
\begin{equation}\label{half surfs}
\begin{aligned}
	\mathcal{S}_{1} &:= \left\{ \Sigma \left( \theta, \phi \right) : \theta \in \left[ \theta_{0} - \pi, \theta_{0} \right), \phi \in \left[ 0, \pi \right] \right\},\quad \textrm{and} \\
	\mathcal{S}_{2} &:= \left\{ \Sigma \left( \theta, \phi \right) : \theta \in \left[ \theta_{0}, \theta_{0} + \pi \right), \phi \in \left[ 0, \pi \right] \right\}.
\end{aligned}
\end{equation}
are symmetric in the sense of section~\ref{section: symmetry}: the symmetric equilibrium to a state in $\mathcal S_1$ is in $\mathcal S_2$ and vice-versa. This same symmetry also implies that $\mathcal S$ has a mirror symmetry about the plane $\phi = \pi/2$.


Because of the mirror symmetry between $\mathcal{S}_{1}$ and $\mathcal{S}_{2}$, the number of relative equilibria for fixed parameters $\Ma$ and $\psi$ is even. We now show that it is generically $0$, $4$, or $8$. Using \eqref{param a} to express $\phi$ in terms of $\Ma$ and $\theta$, and substituting in\re{param psi}, we find
\begin{align*}
	&\left( \frac{\Ma}{\sigma_{1} \sigma_{2}} \left( c_{11} \, \left( \frac{\cos^2 \theta}{\sigma_{1}^2} - \frac{\sin^2 \theta}{\sigma_{2}^2} \right) + c_{12} \, \frac{\cos \theta \, \sin \theta}{\sigma_{1} \, \sigma_{2}}  \right) - \cos \psi \right)^2 \\
	+& \left(\Ma^{2}\left( \frac{\cos^2 \theta}{\sigma_{1}^2} + \frac{\sin ^2 \theta}{\sigma_{2}^2} \right) - 1 \right) \, \left( c_{01} \, \frac{\cos \theta}{\sigma_{1}} + c_{02} \, \frac{\sin \theta}{\sigma_{2}} \right)^2 = 0 \, .
\end{align*}
For fixed values of the parameters $\Ma$ and $\psi$, this is a trigonometric polynomial of degree 4 in $\theta$; it has therefore at most 8 roots \citep[]{Powell1981}.

$\mathcal S_{1}$ does not have boundaries. Indeed, the mapping $(\theta,\phi) \mapsto \Sigma(\theta,\phi)$ is continuous, and therefore any boundary would arise on the borders of the chart $[\theta_{0}-\pi,\theta_{0}) \times [0,\pi]$. However $\mathcal S_{1}$ has no boundary at $\theta \in \{ \theta_{0}-\pi, \theta_{0} \}$: for all $\phi$, we have $\Sigma(\theta_{0} -\pi,\phi) = \Sigma(\theta_{0},\phi)$, and for $\phi$ fixed in $[0,\pi]$, the curve $\{ \Sigma(\theta,\phi): \theta \in [\theta_{0} -\pi,\theta_{0}) \}$ is closed. Furthermore, $\mathcal S_{1}$ has no boundary at $\phi = \{0,\pi\}$: for $\sin \phi = 0$, equation\re{param a} yields $\Ma = 0$, so that for $\phi \in \{0,\pi \}$, the closed curve $\{ \Sigma(\theta,\phi): \theta \in [\theta_{0}-\pi, \theta_{0}) \}$ is reduced to the segment $\{0 \} \times [-\max \{|c_{01}/\sigma_{1}|, \, |c_{02}/\sigma_{2} | \},\max \{|c_{01}/\sigma_{1}|, \, |c_{02}/\sigma_{2} | \}] \times \{ \phi \}$ travelled twice, thereby sealing the boundaries at $\phi \in \{0,\pi\}$. By symmetry, $\mathcal S_{2}$ does not have boundaries either.

Looking for relative equilibria for fixed parameter values $\Ma_{0}$ and $\psi_{0}$ is equivalent to looking for intersections between $\mathcal{S}$ and the straight line in $\left( \Ma,\cos \psi, \phi \right)$-space with $\Ma = \Ma_{0}$ and $\cos \psi = \cos \psi_{0}$. Since $\mathcal S_{1}$ has no boundaries, slicing $\mathcal{S}_{1}$ at $\Ma = \Ma_{0}$ gives a closed curve in the $\left( \cos \psi, \phi \right)$-plane. The straight line $\cos \psi = \cos \psi_{0}$ intersects this curve twice, except at folds or self-intersections of the curve. Since each relative equilibrium is uniquely parametrised by $\left( \theta, \phi \right)$ in the chosen intervals by equations~(\ref{e3 param}, \ref{B param}), self-intersections still count for two distinct relative equilibria. This eventually implies that the parameter space is cut by folds into regions with either 0, 4, or 8 corresponding equilibria (cf. fig.~\ref{subfig: swimmer A mesh top}-\ref{subfig: swimmer A surf top}).

Further geometric properties of surface $\mathcal{S}$ are studied in appendix~\ref{appendix: intersections}.

\subsection{Experimental regimes diagram}
\label{section: exp param}

\begin{figure}[h]
	\centering
	\includegraphics{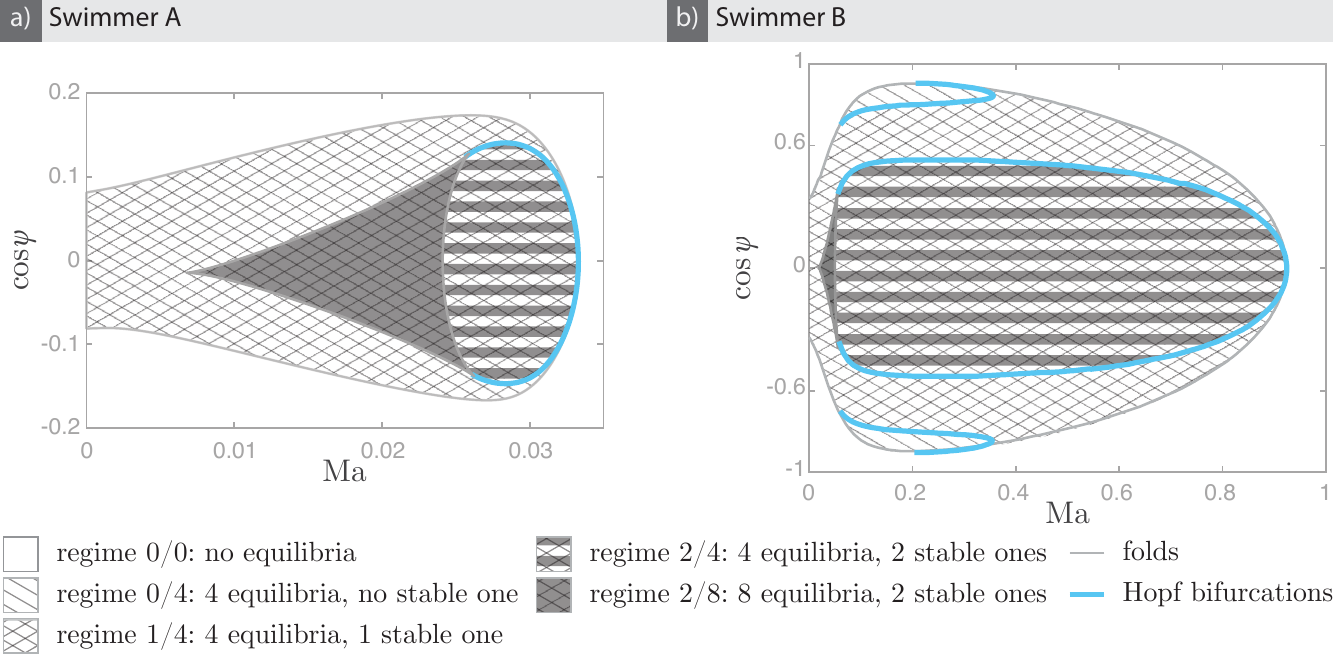}
	{\phantomsubcaption\label{subfig: swimmer A param regimes}}
	{\phantomsubcaption\label{subfig: swimmer B param regimes}}
	\caption[]{Experimental regimes diagrams indicating the number of stable relative equilibria and their total number in function of the two experimental parameters $\Ma$ and $\cos \psi$. Note that the regime 0/4 exists only for swimmer~B. The data shown here corresponds to swimmers~A and~B (cf. section~\ref{section: numerics}).
	}
	\label{fig: param regimes}
\end{figure}

Projecting the surface $\mathcal S$ on the $(\Ma,\cos \psi)$ plane, and keeping track of the number of sheets projected on any given point, we obtain a diagram on which one can read of the number of equilibria (see figures~\ref{subfig: swimmer A mesh top} and~\ref{fig: param regimes}). The general shape of this diagram is consistent across different swimmers: we found at most five different regimes that cover various regions of the $(\Ma, \cos \psi)$ plane. These regions vary in shape and sizes with the shape and magnetisation of a swimmer. However, their general arrangement is consistent as well as their broad features. For instance, the triangular region of the regime 2/8 (standing for 8 equilibria, 2 of which are stable) remains triangular even though it can be of very different proportions for different swimmers. The notable exception is the regime 0/4 present on the wings at large $|\cos \psi|$ for swimmer B and absent for swimmer~A. Some swimmers have it and some don't. When that regime exists, we always found it at similarly large $|\cos\psi|$. 

\section{Helical swimmers and how to operate them}
\label{section: numerics}

\begin{figure}[b!]
	\centering
	\includegraphics{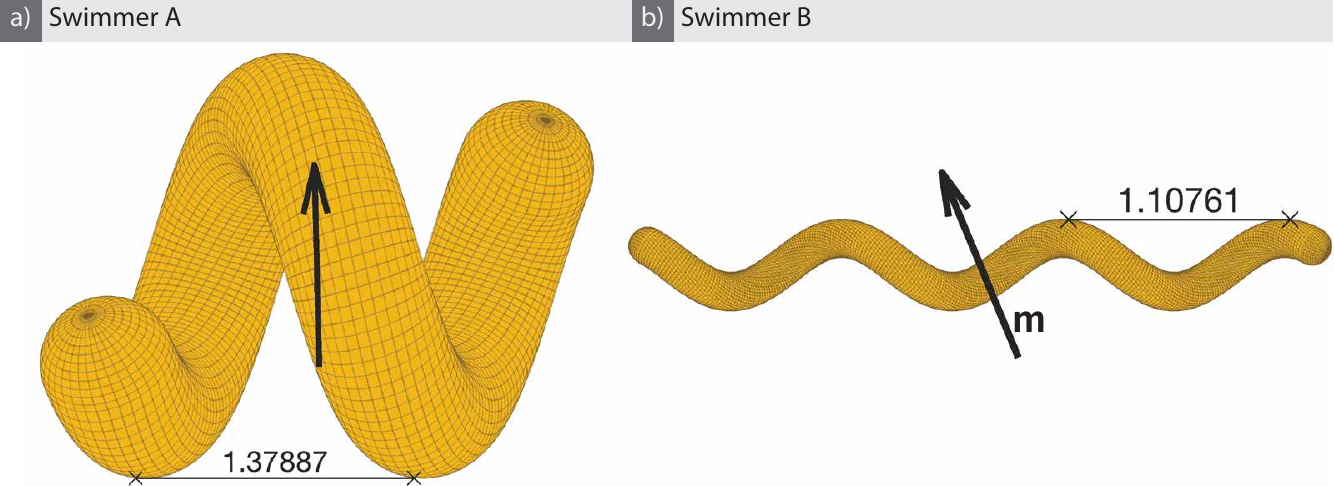}
	{\phantomsubcaption\label{subfig: swimmer A}}
	{\phantomsubcaption\label{subfig: swimmer B}}
	\caption[]{Helical swimmers that we study numerically in detail. The unit magnetic moment $\mbody$ is represented by the arrow, which lies in the paper plane. The detailed numerical data about these swimmers can be found in appendix~\ref{appendix: numerics}.}
	\label{fig: swimmers}
\end{figure}

As stated in the introduction, the principal application motivating our study is the motion of rigid, helical and magnetic microswimmers. We therefore focus on two examples the bodies of which are helical rods with a circular cross-section, capped with half-spheres at the two ends: the swimmers~A and~B shown in fig.~\ref{fig: swimmers}. We analyse many more examples in~\cite{Ruegg-Reymond2019}. To compute the drag matrices $\Drag$ of these helices, we use a code presented in \cite{Gonzalez2009, Li2013, Gonzalez2015} and provided to us by the authors. The outer surface of the helical rod is discretised by one quadrature point in each patch defined by the mesh shown in fig.~\ref{fig: swimmers}. The direction of the magnetic moment $\mbody$ is chosen so that it is not aligned with the helix central axis. Lengths are scaled using the radius of gyration as the characteristic length~$\ell$. The numerical data for $\Drag$ and $\mbody$ corresponding to our examples can be found in appendix~\ref{appendix: numerics}.

\subsection{Experimental regimes diagrams}

The experimental parameters diagrams of both swimmers are shown in figure~\ref{fig: param regimes}. We find 4 regions of parameters exhibiting relative equilibria. It is worth noting that there exist relative equilibria for a much wider range of experimental parameters for swimmer B than for swimmer A, with $\cos \psi$ ranging from $-0.9049$ to $0.9048$ and $\Ma$ up to $0.9244$ compared to $\cos \psi \in \left[ -0.1669, 0.1736 \right]$ and $\Ma \leq 0.0333$ for swimmer A. The different range of mason numbers is readily understood by the analysis of section~\ref{sec-releq}: the maximal Mason admitting a relative equilibrium is given by $\sigma_1$ which is much larger for swimmer~B than for swimmer~A. Note that this is a property of the $\Pbody$ matrix and not of the chirality matrix. We do not have an analytical understanding of the range of $\cos \psi$.

Moreover the regimes diagram~\ref{subfig: swimmer B param regimes} that corresponds to the longer helix is more symmetric about the $\cos \psi = 0$ axis than that of the shorter helix~\ref{subfig: swimmer A param regimes}. It is actually not exactly symmetric. We conjecture that for long helices the diagram will be very close to symmetric regardless of the direction of magnetic moment $\mbody$.

The other major difference between both diagrams is that swimmer~B admits a regime for which relative equilibria exist but are all unstable (regime 0/4). The parameter regions corresponding to this regime are located at the extremal $\cos \psi$ tips of the set of equilibria and are bounded by curves of Hopf bifurcations that have no counterparts for swimmer~A.


Finally, remembering from section~\ref{section: exp param} that these diagrams are obtained by projecting the three-dimensional $\mathcal S$ surfaces, we also include  the corresponding surface from swimmer~B in figure~\ref{fig: swimmer B eq} in the appendix~\ref{appendix: intersections}. Note that for both swimmers, all the stable equilibria are bound to the region $\phi < \frac{\pi}{2}$, in accordance with the intuition that the magnetic moment $\mbody$ tries to align with the magnetic field $\Bbody$. 

Although these diagrams are useful to predict maximal Mason before step-out as well as the number of possible behaviours at given $\Ma$ and $\cos \psi$, they also suggest a number of interesting questions. What happens when the system is multistable? How does the system behave when there are 4 equilibria but none of them are stable? In the simpler case where there exists a single stable equilibrium, are we sure that the system will always collapse to it? How to optimally handle the swimmer to get it where we want it to go? 

To address these questions, we will both apply the techniques developed in section~\ref{section: analysis} and perform numerical experiments.

\subsection{Numerical conditioning and quaternion formulation}
%

We first cast \eqref{ode Q closed} in a form that both allows for easier numerical integration and enables numerical continuations e.g. follow branches of periodic solutions sprouting from Hopf bifurcations as described in section~\ref{section: hopf bifurcations}.

Since the unknown $Q$ in \eqref{ode Q closed} is a rotation matrix, numerical treatment of the system can be eased by using a parametrisation of the group of rotations. Here we parameterise by quaternions which can be treated as 4-vectors by numerical integrators. The quaternion formulation of \eqref{ode Q closed} is \citep[]{Dichmann1996}
\begin{equation}
	\label{quaternion form}
	\begin{aligned}
		\dot{q} &= \frac{1}{2} F^T \left( q \right) \ubody \left( q; \Ma, \psi \right)
		, 
		&
		q \left( 0 \right) &= q_{0},
	\end{aligned}
\end{equation}
where
\begin{equation*}
	F\left( q \right) =
	\begin{bmatrix}
		\phantom{-}q_{4} & -q_{3} & \phantom{-}q_{2} & -q_{1} \\
		\phantom{-}q_{3} & \phantom{-}q_{4} & -q_{1} & -q_{2} \\
		-q_{2} & \phantom{-}q_{1} & \phantom{-}q_{4} & -q_{3}
	\end{bmatrix}
\end{equation*}
and $\ubody \left( q; \Ma, \psi \right) = \Ma \, \ebody_{3} \left( Q\left( q \right) \right) - \omegabody \left( Q\left( q \right); \psi ,\, \Ma\right)$, with
\begin{equation}
	Q \left( q \right) =\left .
	\begin{bmatrix}
		q_{1}^2 -q_{2}^2 -q_{3}^2+q_{4}^2 & 2\left( q_{1}q_{2}-q_{3}q_{4} \right) & 2\left( q_{1}q_{3}+q_{2}q_{4} \right) \\
		2\left( q_{1}q_{2} + q_{3}q_{4} \right) & -q_{1}^2+q_{2}^2-q_{3}^2+q_{4}^2 & 2\left( q_{2} q_{3}-q_{1}q_{4} \right) \\
		2\left( q_{1}q_{3}-q_{2}q_{4} \right) & 2\left( q_{2}q_{3}+q_{1}q_{4} \right) & -q_{1}^2 - q_{2}^2 + q_{3}^2 + q_{4}^2
	\end{bmatrix}\right / |q|^2, \label{paramQ}
\end{equation}
and $\ebody_{3}$, $\omegabody$ given as in \eqref{cl}. This parametersiation of rotations by quaternions is independent of the norm of the quaternion. Accordingly, the Jacobian of the right-hand side of~\eqref{quaternion form} is singular, and this poses a problem for numerical continuation. To circumvent this issue, we use the modified system
\begin{equation}
	\label{ode quat corrected}
	\begin{aligned}
		\dot{q} &= \frac{1}{2} F^T \left( q \right) \ubody \left( q; \Ma, \psi \right) - \frac{1}{2} \left( \left| q \right|^2 - 1 \right) q
		, &
		q\left( 0 \right) &= q_{0}.
	\end{aligned}
\end{equation}
 This modified system ensures that the solutions are unit quaternions and the stability is the same as for the original system, both for steady states and periodic solutions.

\subsection{Periodic solutions}
\label{section: periodic}

\begin{figure}[b!]
	\centering
	\includegraphics{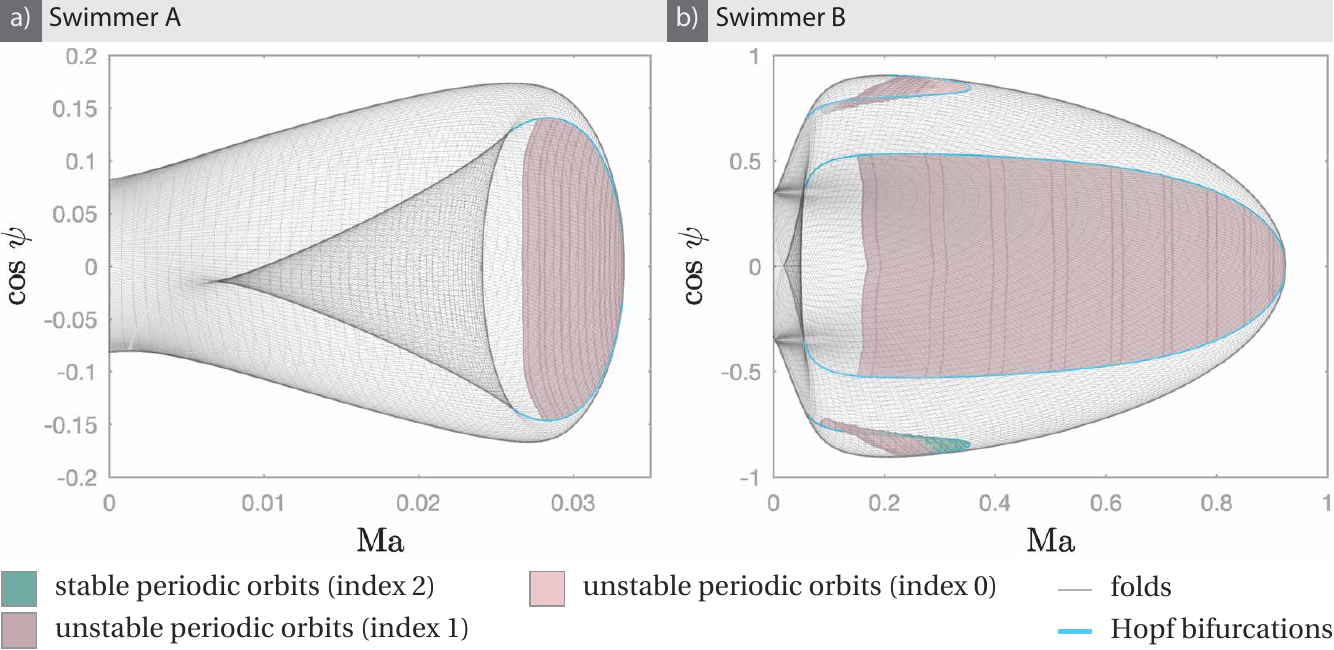}
	{\phantomsubcaption\label{subfig: swimmer A mesh periodic}}
	{\phantomsubcaption\label{subfig: swimmer B mesh periodic}}
	\caption[]{Foils of periodic orbits branching from Hopf bifurcations on the half-surfaces $\mathcal S_{1}$ shown in figures~\ref{subfig: swimmer A mesh top} and~\ref{subfig: swimmer B mesh top} respectively. Only half $\mathcal S_{1}$ is shown; a similar foil is obtained from the symmetric curve of Hopf bifurcations on $\mathcal S_{2}$. The stability is reversed on the symmetric half $\mathcal S_{2}$ not shown. The data shown here corresponds to the helical swimmers~A and~B (cf. fig.~\ref{fig: swimmers}).
	}
	\label{fig: periodic}
\end{figure}

Because the system displays Hopf bifurcation, we expect that it has periodic orbits for some regions of the experimental parameters $\Ma$ and $\cos \psi$. We tracked these solutions with the Matlab package MatCont~\citep[]{Dhooge2008} applied to the ODE \eqref{ode quat corrected}. The stability index of these periodic orbits is obtained from the Floquet multipliers~\citep[]{Kuznetsov2004} which are part of the code output. The starting points for continuations are chosen on the curves of Hopf bifurcations (drawn in blue in all figures). Each point on those curves is associated with a period depending on the purely imaginary eigenvalue of the stability matrix at that point. Then we chose to continue by finding the $\Ma$ and $\cos \psi$ that keep the period of the periodic orbit constant. By this process, each point on the curves of Hopf bifurcation gives rise to a curve of periodic orbits in the $(\Ma,\cos\psi)$ plane. We found that these curves both start and end at Hopf bifurcations. Gathering these curves of periodic orbits into a set forms a sheet of periodic orbits as shown in red (if unstable) and green (if stable) in figure~\ref{fig: periodic}. These sheets are obtained by interpolation of the darker lines, which represent branches obtained by numerical continuation of system\re{ode quat corrected} with period kept constant.

Note that at the beginning of this process, each point along a curve of Hopf bifurcations is associated with a period. In our system, all curves of Hopf bifurcations were found to begin and end at folds. This implies that the associated period goes to infinity making accurate numerical continuations difficult in the vicinity of these points. This is why the foils shown in figure~\ref{fig: periodic} are incomplete.

For both swimmers A and B, the foil originating from the Hopf curve bordering the 2/4 regime seems to fill the very domain of that regime (we could not have a definite answer because of the numerical issue mentioned above). In our examples, all periodic solutions existing for regime 2/4 are unstable.

For swimmer B, the foil originating from the Hopf curve bordering the 0/4 regime seems to fill the domain of that regime but it is not limited to it. We note that some of these periodic orbits are stable and some are unstable. It is noteworthy that this foil is not restricted to the 0/4 regime as some of the unstable branches extend into the 1/4 regime.

There also exist other types of stable periodic orbits that do not connect to Hopf bifurcations. In particular, we study stable periodic orbits in regions which are well away in the 0/0 region and in the $\Ma\ll1$, $\Ma \gg 1$ and $\sin \psi\ll 1$ limits in a separate paper~\citep{Ruegg-Reymond2018}. These results are also described in the thesis of \cite{Ruegg-Reymond2019}.

In conclusion, as long as we stay in the regimes where there exist stable equilibria and well away from large $|\cos\psi|$, we systematically observed that all periodic orbits were unstable.

\begin{figure}[p!]
	\centering
	\includegraphics{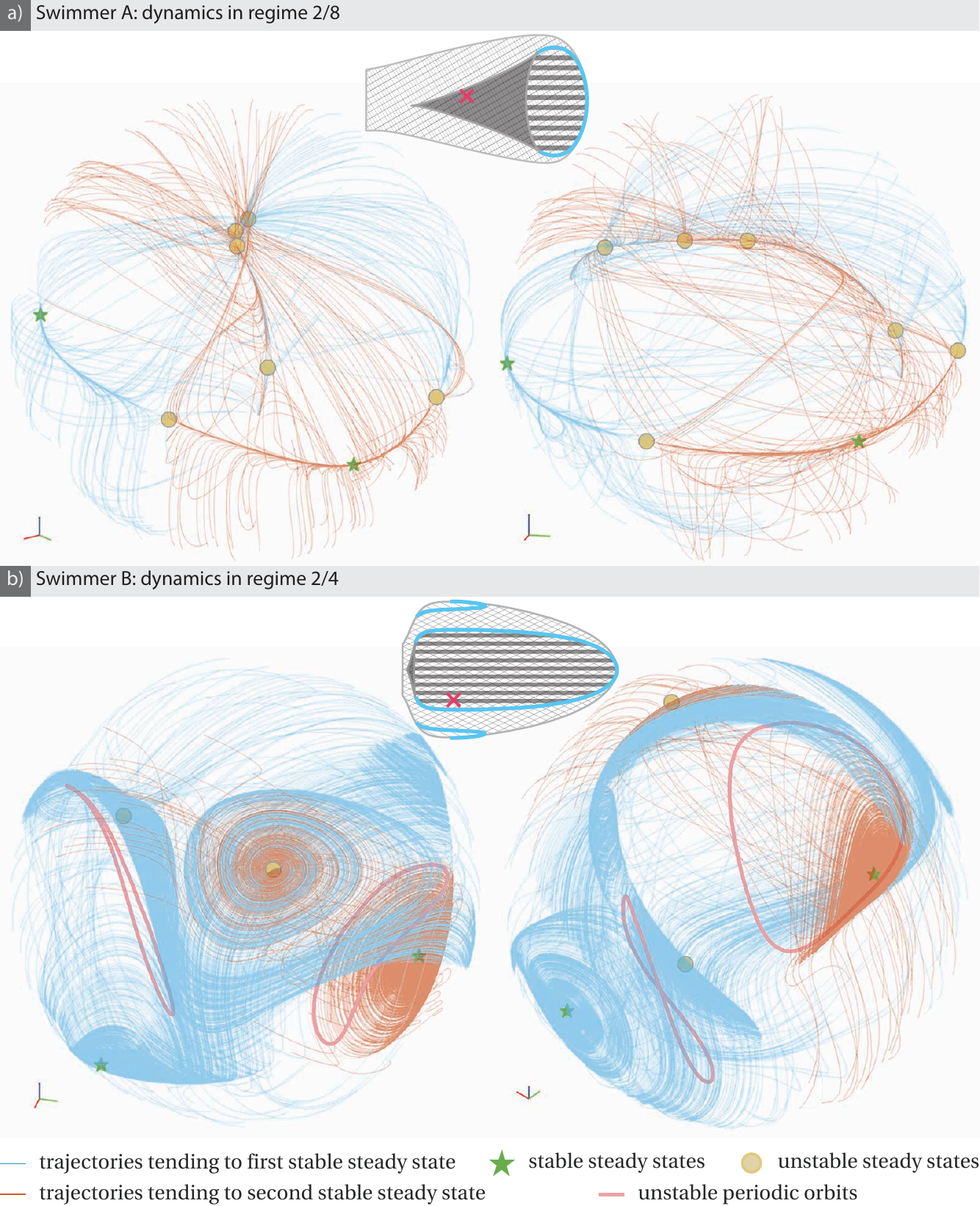}
	{\phantomsubcaption\label{subfig: dynamics A 2/8}}
	{\phantomsubcaption\label{subfig: dynamics B 2/4}}
	\caption[]{ Phase portrait of the dynamical system \eqref{quaternion form} describing the orientation of a swimmer in time by means of quaternions. Unit quaternions are represented by their vector part lying in the unit ball. The data shown here corresponds to \subref{subfig: dynamics A 2/8}) swimmer~A with parameter values $\Ma = 0.015$ and $\cos \psi = 0.01$, and \subref{subfig: dynamics B 2/4}) swimmer B with $\Ma = 0.2198$ and $\cos\psi = -0.3989$.  The chosen parameter values are marked by a cross on the corresponding experimental regimes diagrams. On each row the same object is represented under two different viewpoints.
	}
	\label{fig: dynamics}
\end{figure}

\subsection{Phase portrait and basins of attraction}
\label{section: behaviour at specific parameters}

For both swimmers A and B, we study phase portraits of the dynamical system \eqref{quaternion form}  obtained by integrating it numerically with many initial conditions for various parameter values $(\Ma,\, \psi)$. Showing quaternion dynamic is not geometrically straightforward as we have to project a four-dimensional dynamic onto a 2D sheet. The matter is simplified by the fact that the parameterisation~\ref{paramQ} is independent of the norm of the quaternion: we can therefore drop one degree of freedom. The unit quaternions are represented in 3D by their three vector components $q_{1}$, $q_{2}$, and $q_{3}$, and since opposed quaternions represent the same rotation, we have the convention that the scalar part $q_{4} \ge 0$. This representation is the 4D equivalent of representing the upper hemisphere of the sphere $S^2$ by flattening it onto the unit disc. Due to the equivalence of opposite quaternions, antipodal points on the outer boundary are equivalent. The phase portraits displayed in figure~\ref{fig: dynamics} show the vector part of the quaternions so that the dynamic is actually a dynamic in the unit ball.

Regime 1/4 is the only observed regime for which there is a unique stable steady state. We observed two different behaviours with numerical experiments in this regime. In most cases, there is no other stable special solution, and all trajectories obtained by numerical integration of\re{ode quat corrected} reach the stable steady state after some time, suggesting that it is globally attractive. For some swimmers, a stable periodic orbit coexists with a stable steady state for parameter regime 1/4. Each stable solution then has its own attraction basin. This second behaviour is studied in more detail in~\citep{Ruegg-Reymond2019}.

In figure~\ref{fig: dynamics}, we look closely at two particular examples of phase portraits for parameter values corresponding to cases where there are two stable relative equilibria: we picked a 2/8 regime for swimmer~A ($\Ma = 0.015$ and $\cos \psi = 0.01$) and a 2/4 regime for swimmer~B ($\Ma = 0.2198$ and $\cos\psi = -0.3989$). In both cases, all trajectories were found to tend to one or the other stable relative equilibrium which indicates that the union of their basins of attraction is likely to cover almost all of the possible initial orientations. Note that the trajectories are coloured according to the equilibrium that is reached so the blue and red shapes effectively draw out these basins of attraction. In regime 2/8 for swimmer~A (fig.~\ref{subfig: dynamics A 2/8}) we can observe that the two relative equilibria lie on two different accumulations of lines prefiguring the existence of two heteroclinic orbits, each linking together four relative equilibria. In regime 2/4 for swimmer~B (fig.~\ref{subfig: dynamics B 2/4}) the four relative equilibria coexist with two unstable periodic orbits found by numerical continuation from Hopf bifurcations (cf. section~\ref{section: periodic}).

\subsection{Optimising the magnetisation}

\begin{figure}[t!]
	\includegraphics{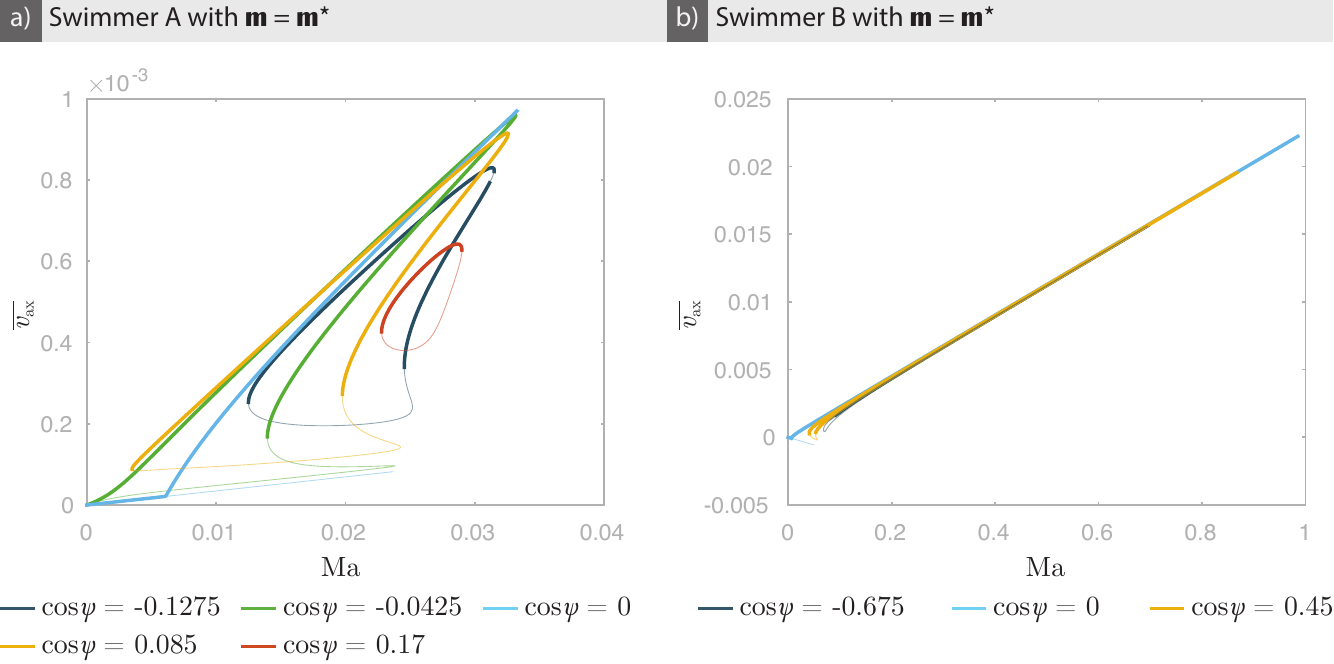}
	{\phantomsubcaption\label{subfig: axvelVsMa a opt}}
	{\phantomsubcaption\label{subfig: axvelVsMa b opt}}
	\caption{Axial Velocity vs $\Ma$ for several values of the conical angle $\psi$ for swimmers with optimal magnetisation direction (cf. table~\ref{table: opt axvel}).  Thinner parts of the curves correspond to values reached only for unstable relative equilibria. The data shown here corresponds to swimmers with the same shape as swimmer A and  swimmer B (cf. fig.~\ref{fig: swimmers}).}
	\label{fig: axvelVsMa opt m}
\end{figure}

Following the procedure described in section~\ref{section: opt m}, we look for a magnetisation direction $\mbody$ that maximises the axial velocity $\vax$ over relative equilibria for swimmers with the same shapes as the helical swimmers A and B (cf.~fig.~\ref{fig: swimmers}). Results are summarised in table~\ref{table: opt axvel}. The optimal magnetic moment direction $\mbody^\star$ is given with respect to the basis described in appendix~\ref{appendix: numerics} -- the third director is aligned with the helix axis. Axial velocities as functions of $\Ma$ for various values of $\cos\psi$ of the optimally magnetised versions of swimmers A and B are shown in figure~\ref{fig: axvelVsMa opt m}.

\begin{table}[h!]
	\centering
	\footnotesize{
	\begin{tabularx}{\textwidth}{l C C C C}
		\toprule
		Swimmer geometry & $\mbody^\star$ & $\vax^\star$ & \multicolumn{2}{c}{Experimental parameters} \\
		\cmidrule(l){4-5}
		&&& $\Ma^\star$ & $\psi^\star$ (rad) \\
		\hline
		Similar to swimmer A & $\begin{bsmallmatrix}-0.9833 \\ \phantom{-}0.0000 \\ -0.1819\end{bsmallmatrix}$ & $9.7261 $ $\times 10^{-4} \, \frac{\ell}{t_{c}}$ & 0.0333 & 1.5708 \\
		Similar to swimmer B & $\begin{bsmallmatrix}\phantom{-}0.0000\\\phantom{-}0.9955\\-0.0949\end{bsmallmatrix}$ & $0.0223 \, \frac{\ell}{t_{c}}$ & 0.9880 & 1.5708 \\
		\bottomrule
	\end{tabularx}}
	\caption{Optimal magnetisation direction $\mbody^\star$ for two swimmer geometries, with associated optimal axial velocity and experimental parameters.}
	\label{table: opt axvel}
\end{table}

\subsection{How to operate bistable swimmers?}

We observe that for arbitrary magnetisations, optimal experimental parameters are typically found in the bistable regime. Furthermore, we provided examples in section~\ref{section: behaviour at specific parameters} for which the basins of attraction of the two stable relative equilibria (almost) fill the full phase space. In practice, it would therefore not do to simply switch on the external field at the optimal $\Ma$ and $\cos \psi$. Indeed by doing so, one would sometimes reach the optimal regime and sometimes reach the other stable equilibrium. Here we discuss how to systematically bring the swimmer in one chosen relative equilibrium.

\begin{figure}[t!]
	\centering
	\includegraphics{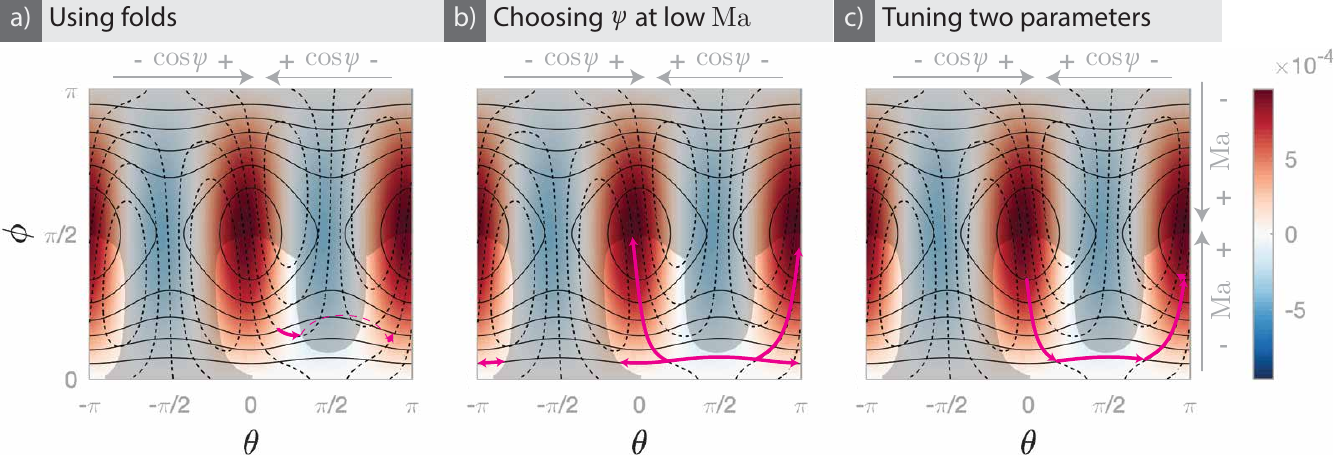}
	{\phantomsubcaption\label{subfig: axvel path fold}}
	{\phantomsubcaption\label{subfig: axvel path low Ma}}
	{\phantomsubcaption\label{subfig: axvel path 2 params}}
	\caption{Different ways to recover the other stable state in bistable modes: \subref{subfig: axvel path fold}) decrease $\cos \psi$ until a fold point and increase it again once the other branch is reached; \subref{subfig: axvel path low Ma}) choose $\cos \psi$ with a low value of $\Ma$ so as to reach the desired stable branch once $\Ma$ is increased; \subref{subfig: axvel path 2 params}) Decrease $\Ma$, increase $\cos \psi$, increase $\Ma$, and then decrease $\cos \psi$. The data used here corresponds to swimmer~A (cf. fig.~\ref{fig: swimmers})}
	\label{fig: axvel paths}
\end{figure}

We present two strategies shown  in figure~\ref{fig: axvel paths}. In all cases, they should work without requiring detailed knowledge of the particular swimmer. The first strategy is to use the fact that the region of stability of the two equilibria are in most cases different. Hence, a simple try is to keep $\Ma$ constant and to change $\psi$ slowly until the current branch is destabilised. Then, we can slowly bring $\psi$ back to its original value: if no further jump is observed then the swimmer is now at the (only) other equilibrium -- see the panel~\subref{subfig: axvel path fold}) in figure~\ref{fig: axvel paths}. 

The other strategy is based on the following observations. First, in all examples presented here and in \cite{Ruegg-Reymond2019}, we found that the stable region of the $(\theta,\phi)$ plane is simply connected (when accounting for the $2\pi$-periodicity in $\theta$). This means that it must be possible to smoothly transform any stable equilibrium in another equilibrium by smoothly changing $\theta$ and $\phi$ or equivalently by smoothly changing $\Ma$ and $\psi$. Next, direct inspection of~\eqref{param a psi} shows that one can qualitatively understand the transition from one to two stable equilibria by considering the low $\sin \phi$ and $\sin\phi\to 1$ limits (keeping in mind that stable equilibria are found only for $\phi\leq\pi/2$. In the first case, $\Ma$ is small and the first term dominates in~\eqref{param psi}: there is a two-to-one relation from $\theta$ to $\psi$. But in the second case, $\Ma$ is large and the second term dominates in~\eqref{param psi}: there is a four-to-one relation from $\theta$ to $\psi$. In practice, we have systematically observed that only one of the two equilibria at $\sin\phi\ll 1$ is ever stable and that two of the equilibria at $\sin \phi\to 1$ are stable.

Therefore, at low $\Ma$, we have low\footnote{That is because the second factor in~\eqref{param a} is bounded from below by $\sigma_2$ (the smaller non-vanishing singular value of $\Pbody$). Hence $\Ma \ll \sigma_2 \Longrightarrow\sin\phi\ll1$.} $\sin \phi$ and in that region, there is a one-to-one relation between $\cos\psi\in I =  [-\sqrt{1-(\betabody_0.\etabody_0)^2},\sqrt{1-(\betabody_0.\etabody_0)^2}] $ and the unique stable equilibria existing at that value of\footnote{Note that the range $I$ can be very narrow for some swimmers and fairly broad for others. For example, for swimmer $A$ we have $\sqrt{1-(\betabody_0.\etabody_0)^2}\simeq 0.08$ while for swimmer B we have $\sqrt{1-(\betabody_0.\etabody_0)^2}\simeq0.3$. That is for swimmer A, equilibria will be found at low $\Ma$ provided that $\psi$ deviates from $90^\circ$ by no more than $\sim 5^\circ$. While for swimmer B, $\psi$ can deviate  from $90^\circ$ by up to $\sim 17^\circ$.} $\Ma$.

The second strategy is therefore as follows (see also panels~\subref{subfig: axvel path low Ma}-\subref{subfig: axvel path 2 params}) of figure~\ref{fig: axvel paths}): 
\begin{enumerate}[label=Step \arabic*\,:\,, itemindent=38pt+\parindent,itemsep=0pt,leftmargin=0pt,topsep=0pt]
\item Study the swimmer at low $\Ma$. Starting at $\cos\psi =0$ which is always stable, slowly increase $\psi$ up to the loss of stability: this defines the boundary of the $I$ interval above.
\item Bring $\cos \psi$ well on one side of the $I$ interval (there is no need to be on the boundary).
\item  Slowly increase $\Ma$, then optimise $\Ma$ and $\cos \psi$ by trial and errors. Note the optimal values $\Ma^\star$ and $\cos\psi^\star$. At this stage there is no guarantee that we are at the optimal equilibrium: the other stable one might be better.
\item Slowly bring $\cos \psi$ back to zero and then slowly decrease $\Ma$ back to small values.
\item Set $\cos \psi$ to the other  side of $I$ (as compared to the choice made at step 1). 
\item Slowly increase $\Ma$ back to $\Ma^\star$, then slowly bring $\cos \psi$ back to $\cos \psi^\star$ and compare the axial velocity obtained in this final state to that which was obtained at step 3.
\item If the last axial velocity is greatest then locally optimise $\Ma$ and $\cos\psi$ again and you found the optimum. If the axial velocity of step 4 was higher, then come back to that state by following backward the sequence 4 to 6.
\end{enumerate}

\section{Axial velocity of the three beads cluster}
\label{sec: comparison}

\begin{figure}[b!]
	\includegraphics{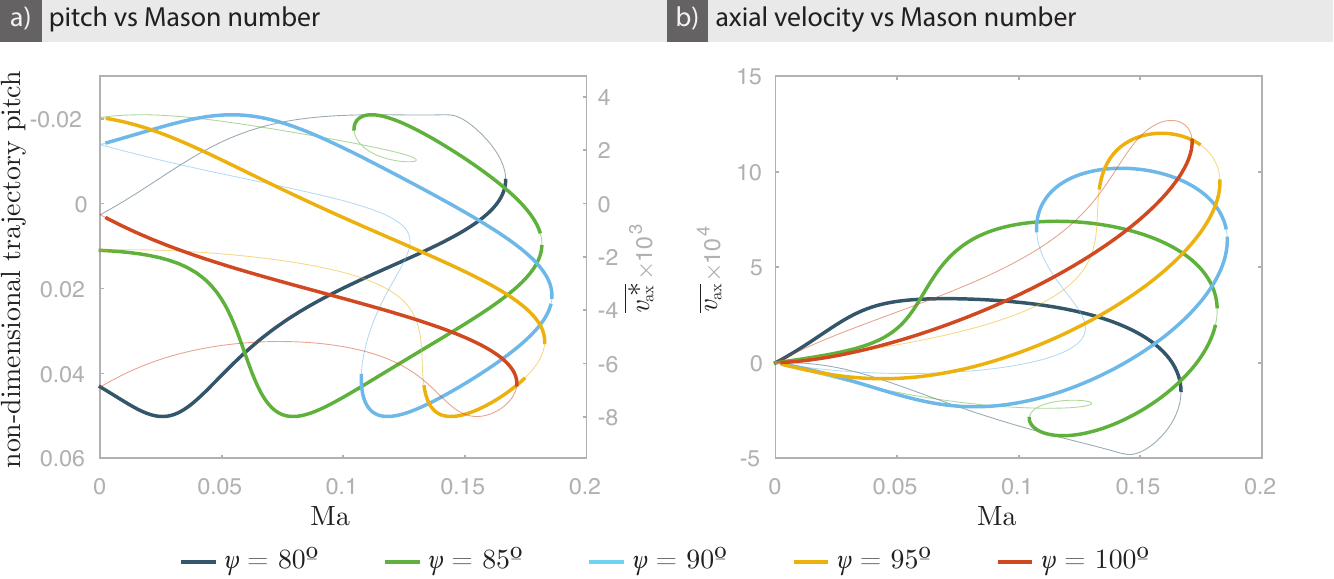}
	{\phantomsubcaption\label{subfig: meshkati pitch}}
	{\phantomsubcaption\label{subfig: meshkati axvel}}
	\caption[]{Trajectory pitches and axial velocities plotted against Mason number for different values of $\psi$ for the three beads cluster from \cite{Meshkati2014}. The thin lines correspond to unstable relative equilibria. The two scales indicated on the vertical axis of panel~\subref{subfig: meshkati pitch}) correspond to two different scalings: the scaling that we introduced in section~\ref{section: dim analysis} on the left, and the scaling presented by~\cite{Meshkati2014} on the right. This figure can directly be compared with figure~11 of \cite{Meshkati2014}.}
	\label{fig: meshkati comp}
\end{figure}
\begin{figure}[t!]
	\includegraphics{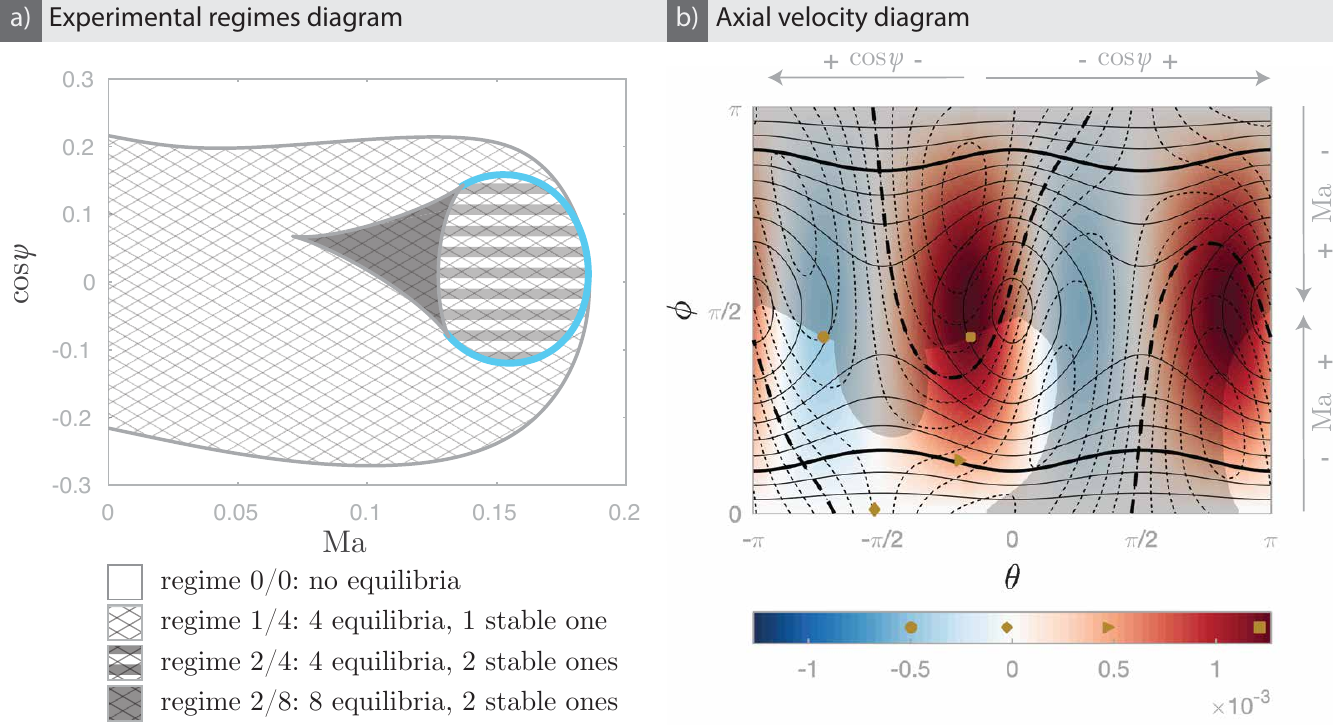}
	{\phantomsubcaption\label{subfig: meshkati regimes}}
	{\phantomsubcaption\label{subfig: meshkati contours}}
	\caption[]{Three beads cluster from \cite{Meshkati2014}: \subref{subfig: meshkati regimes}) Regime diagram in function of experimental parameters $\Ma$ and $\cos \psi$. \subref{subfig: meshkati contours}) Axial velocities of all relative equilibria in terms of the mapping parameters $\theta$ and $\phi$. Level curves of $\Ma$ are indicated in solid lines (blod line: $\Ma = 0.06$), and level curves of $\psi$ in dashed lines (bold line: $\cos \psi = -0.05$).}
	\label{fig: meshkati regimes+contours}
\end{figure}

Many studies report experimentally and/or numerically on the dependance of the axial velocity (up to scaling choices) on the Mason number (up to scaling choices) at constant angle $\psi$ between the external field and its axis of rotation. The parameterisations of the axial velocity~\eqref{vax} and of $\Ma$ and $\cos \psi$~\eqref{param a psi} proposed here provide an analytical expression of that relation. Indeed, the curve $\vax(\Ma)$ is recovered by looking up the values of $\Ma$ and $\vax$ along a level curve of the function $\cos\psi(\theta,\phi)$. In this section, we use this new tool to revisit swimmers that were used as examples in \cite{Meshkati2014} and in~\cite{Morozov2017}. It consists in three contiguous spherical beads arranged so that the lines joining their centres form a $90^\circ$ angle. 
 
\begin{figure}[b!]
	\includegraphics{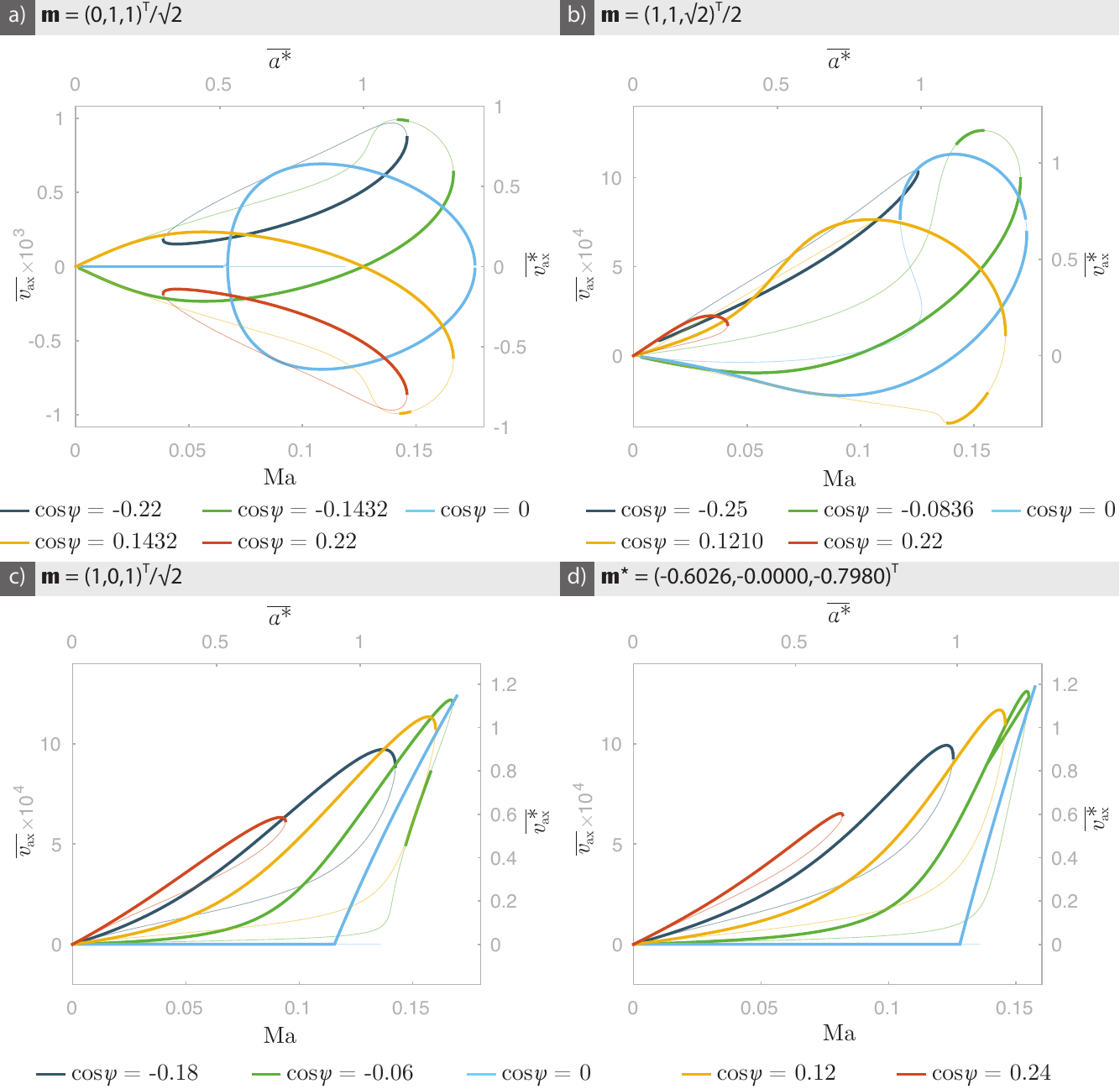}
	{\phantomsubcaption\label{subfig: morozov 3b m1}}
	{\phantomsubcaption\label{subfig: morozov 3b m3}}
	{\phantomsubcaption\label{subfig: morozov 3b m2}}
	{\phantomsubcaption\label{subfig: morozov 3b mstar}}
	\caption[]{Three-bead cluster of \cite{Morozov2017} with a 90$^\circ$ vertex angle, for several magnetic moment directions $\mbody$: axial velocities plotted against Mason number for different values of $\psi$. Thin lines correspond to unstable relative equilibria. Two scales are indicated: $\Ma$ and $\overline \vax$ refer to the Mason number and non-dimensional axial velocity as introduced in section~\ref{section: dim analysis}, while $\overline{\alpha^*}$ and $\overline{v_\text{ax}^*}$ refer to the non-dimensional angular velocity of the magnetic field and non-dimensional axial velocity using the scaling proposed by~\cite{Morozov2017}. Panels~\subref{subfig: morozov 3b m1}) and~\subref{subfig: morozov 3b m2}) can be compared with figure~7 of \cite{Morozov2017}, and panel~\subref{subfig: morozov 3b m3} with their figure~9.}
	\label{fig: morozov 3b comp}
\end{figure}

To compute, we rescale the mobility matrices provided in the original studies as described section~\ref{section: dim analysis}. The resulting non-dimensional matrices are listed in appendix~\ref{appendix: numerics}. Our results are gathered in figures~\ref{fig: meshkati comp} and~\ref{fig: morozov 3b comp}, where we systematically provide two sets of axes: all the starred quantities correspond to the scaling choices of the original authors while the plain ones are as defined in the present document. A bold line indicates the presence of a branch of stable equilibria, while a thin line shows unstable equilibria. To understand the number of equilibria in the system, it is important to note that the axial velocity of an equilibrium and of its symmetric twin (in the sense of section~\ref{section: symmetry}) are equal (see equation~\eqref{ax vel} together with the fact that at equilibria, $\breve \elab_3 = -\elab_3$). Hence every point on these plots correspond to a \emph{pair} of equilibria. As already mentioned at most one equilibrium of a symmetric pair can be stable.

\begin{figure}[b!]
	\includegraphics{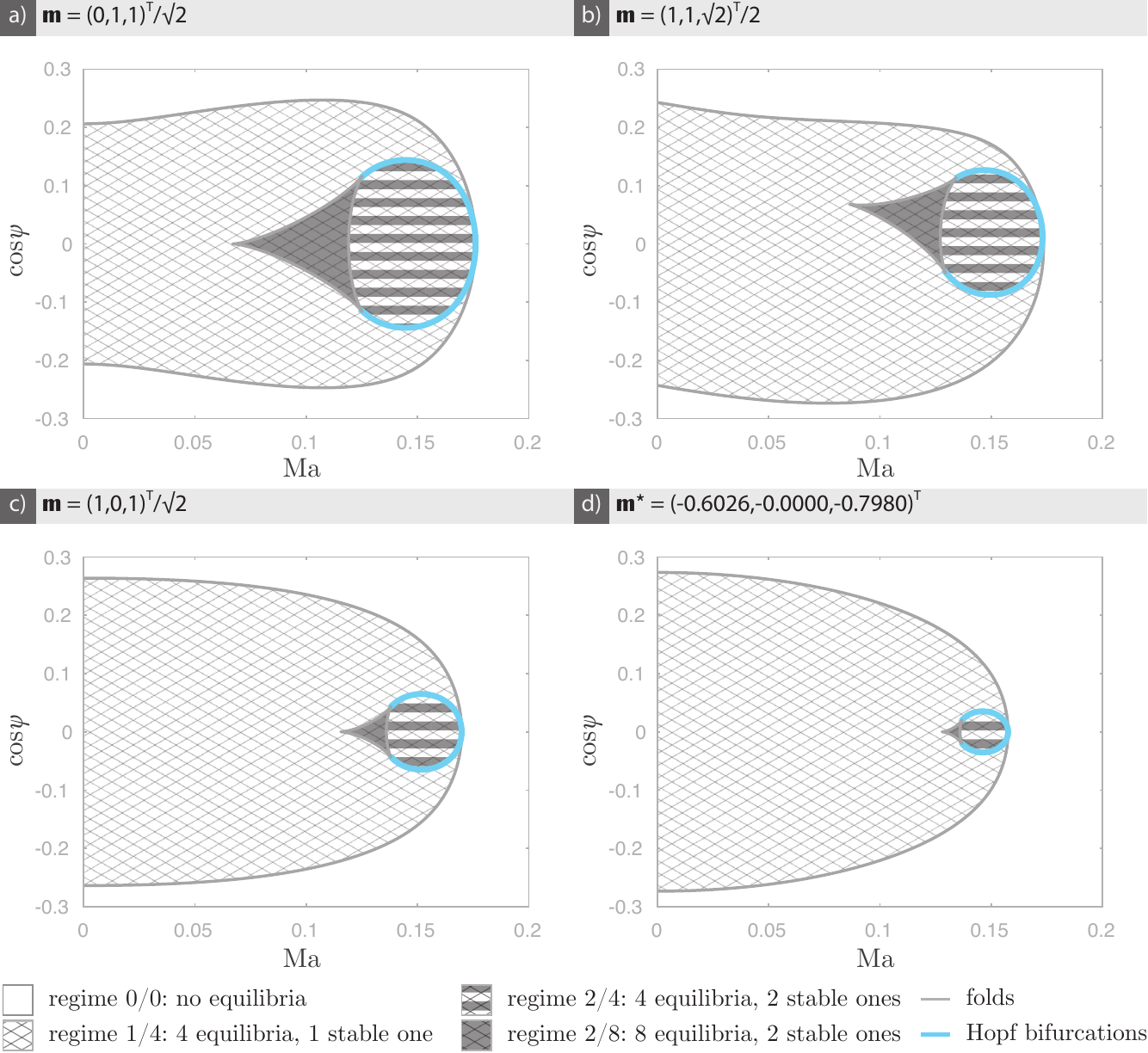}
	{\phantomsubcaption\label{subfig: morozov 3b m1 regimes}}
	{\phantomsubcaption\label{subfig: morozov 3b m3 regimes}}
	{\phantomsubcaption\label{subfig: morozov 3b m2 regimes}}
	{\phantomsubcaption\label{subfig: morozov 3b mstar regimes}}
	\caption[]{Three-bead cluster of \cite{Morozov2017} with a 90$^\circ$ vertex angle, for several magnetic moment directions $\mbody$: diagram of the experimental parameters regimes.}
	\label{fig: morozov 3b regimes}
\end{figure}

\begin{figure}[b!]
	\includegraphics{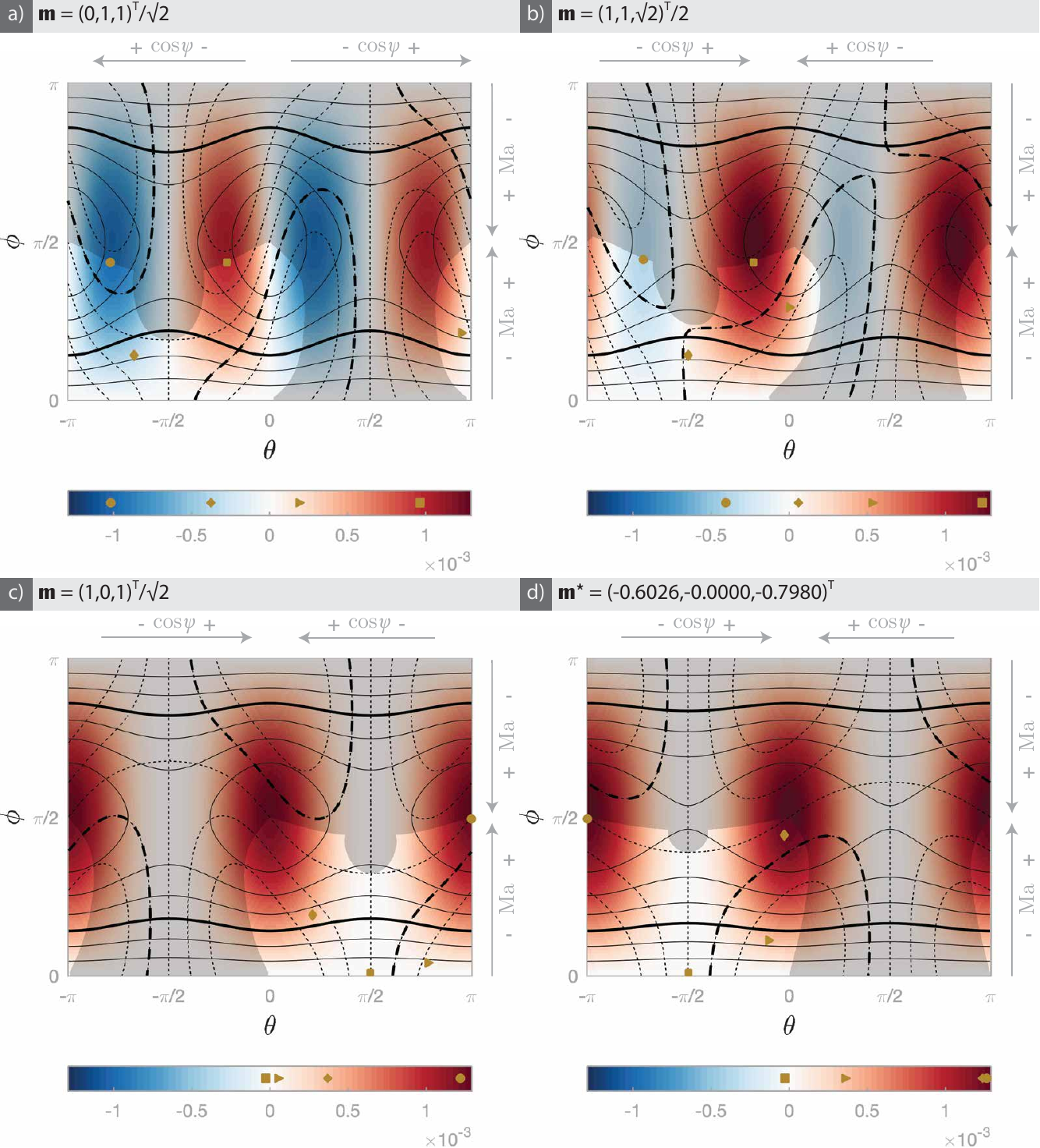}
	{\phantomsubcaption\label{subfig: morozov 3b m1 contours}}
	{\phantomsubcaption\label{subfig: morozov 3b m3 contours}}
	{\phantomsubcaption\label{subfig: morozov 3b m2 contours}}
	{\phantomsubcaption\label{subfig: morozov 3b mstar contours}}
	\caption[]{Three-bead cluster of \cite{Morozov2017} with a 90$^\circ$ vertex angle, for several magnetic moment directions $\mbody$: axial velocity as a function of the mapping parameters $\theta$ and $\phi$ with level curves of the experimental parameters $\Ma$ in solid lines and $\cos\psi$ in dashed lines. Highlighted contours are \subref{subfig: morozov 3b m1 contours}) $\Ma = 0.0754$ and $\cos \psi = 0.0823$, \subref{subfig: morozov 3b m3 contours}) $\Ma = 0.0743$ and $\cos \psi = 0.0708$, \subref{subfig: morozov 3b m2 contours}) $\Ma = 0.0728$ and $\cos \psi = 0.0880$, \subref{subfig: morozov 3b mstar contours}) $\Ma = 0.0675$ and $\cos \psi = 0.0912$.}
	\label{fig: morozov 3b contours}
\end{figure}

Figures~\ref{subfig: meshkati regimes} and~\ref{fig: morozov 3b regimes} show the corresponding experimental regime diagrams from which one can readily read the number of (stable) equilibria for any pairs of experimental parameters values of $\Ma$ and $\cos \psi$, while figures~\ref{subfig: meshkati contours} and~\ref{fig: morozov 3b contours} display the $(\theta,\phi)$-plane with axial velocity and level curves of $\cos\psi$ and $\Ma$.

Generally, \cite{Morozov2017} studied the axial velocity as a function of the shape and magnetisation of the swimmer. In particular, they showed that for the same geometry, different magnetisation orientation can produce vastly different responses: for some, increasing the Mason number provokes a linear increase of the swimmer's velocity (e.g. the curve $\cos\psi = 0$ replicated on fig.~\ref{subfig: morozov 3b m2}), while for others, the response is not linear (e.g. the curve $\cos \psi = 0$ replicated on fig.~\ref{subfig: morozov 3b m1}, see also~\cite{Fu2015}). The authors also discussed the optimisation of both the shape and the magnetisation in some specific cases. However their analysis was focused on the case of $\psi = \pi/2$. It turns out that faster axial motion can be achieved by relaxing the constraint.


\cite{Meshkati2014} are to our knowledge the only other authors studying cases with $\psi\neq \pi/2$. They obtained an algebraic system of equations for the relative equilibria equivalent to our~\eqref{equilibrium conditions} and proceeded to solve it numerically. Interestingly, they reported gaps in the $\vax(\Ma)$ curves which they attributed to numerical issues with their root finding algorithm. The present study indicates that the branch of equilibria does exist but that these equilibria are actually unstable precisely in the gaps where their algorithm failed. Hence the axial velocities in those gaps are not achievable in actual experiments.

The axial velocities obtained with the mobility matrix computed by \cite{Meshkati2014,Meshkati2017} are shown in figure~\ref{fig: meshkati comp} in a way that can be straightforwardly compared with figure~11 of \cite{Meshkati2014}. We keep their characteristic length $\ell$, chosen as the distance between the two furthest bead centres, expressed in terms of the radius of one bead $R_{b}$ as $\ell = 2 \sqrt{2} R_{b}$. Because the authors used the inverse axial velocity $1/\alpha$ of the magnetic field as a characteristic timescale, their non-dimensional axial velocity $\overline{\vax^*}$ is proportional to our pitch $\overline{p}$. There is also a sign difference between their axial velocity $\overline{\vax^*}$ and our pitch $\overline{p}$. It might be due to different sign conventions, but we have not been able to pinpoint the exact cause.  Both scales are shown on fig.~\ref{subfig: meshkati pitch} to make comparison with their figure~11 easier. Note that the dimensional axial velocity $\vax$ of the swimmer can be recovered from our non-dimensional axial velocity $\overline \vax$ as $\vax = \overline \vax \, m \, B/(\eta \, \ell^2)$, and from their non-dimensional axial velocity $\overline{\vax^*}$ as $\vax = \overline{\vax^*} \, \ell \, \alpha$, while the dimensional trajectory pitch $p$ is recovered from the non-dimensional trajectory pitch $\overline p$ as $p = \ell \, \overline p$.


\cite{Morozov2017} report using the same algorithm as \cite{Meshkati2014} to compute the mobility matrix of the three bead swimmer, without explicitly presenting their result. To compare our results, we transformed the mobility matrix given by \cite{Meshkati2014} according to the conventions of \cite{Morozov2017} (cf. appendix~\ref{appendix: numerics}. We report the corresponding axial velocities in figure~\ref{fig: morozov 3b comp} for a few choices of magnetisation. This figure is to be compared with figures 7 and 9 of \cite{Morozov2017}. To make comparison straightforward, we indicated both the scales we use, and the scales they chose. The dimensional axial velocity $\vax$ can be recovered from our non-dimensional velocity $\overline \vax$ as $\vax = \overline \vax \, m \, B/(\eta \, \ell^2)$, where $\ell = 2\sqrt{2} R_{b}$, with $R_{b}$ the radius of one bead, or from their non-dimensional axial velocity $\overline{\vax^*}$ as $\vax = \overline{\vax^*} \, m \, B \, |\mathrm{Ch}| \, F_{\bot}/(\eta \, R_{b}^2)$, where $|\mathrm{Ch}|$ and $F_{\perp}$ are defined as in appendix~\ref{appendix: numerics}. The dimensional angular velocity of the magnetic field $\alpha$ is recovered from their non-dimensional axial velocity as $\alpha = \overline{\alpha^*} \, m \, B \, F_{\bot}/(\eta \, R_{b}^3)$ -- note it can also be computed from the Mason number using formula~\ref{Mason number}.

It is well worth noting that for a given magnetisation, the optimal axial velocities are often obtained for values of $\psi$ well away from $\pi/2$. In particular, for $\mbody = (0,1,1)^T/\sqrt{2}$ (fig.~\ref{subfig: morozov 3b m1}), the extremal non-dimensional axial velocities
, which are obtained for $\cos \psi = \mp 0.1432$ and $\Ma = 0.1433$, are 26\% larger (in absolute value) than the extremal axial velocities reached with the constraint $\cos \psi = 0$. Similarly, for $\mbody = (1,1,\sqrt{2})^T/2$ (fig.~\ref{subfig: morozov 3b m3}), the maximal non-dimensional axial velocity, reached with $\cos \psi = -0.0836$ and $\Ma = 0.1525$, is 13\% larger than the maximal axial velocity reached with $\cos \psi = 0$. The last panel of figure~\ref{fig: morozov 3b comp} showcases the optimal magnetisation obtained by our method of section~\ref{section: opt m}.

\begin{figure}[p!]
	\includegraphics{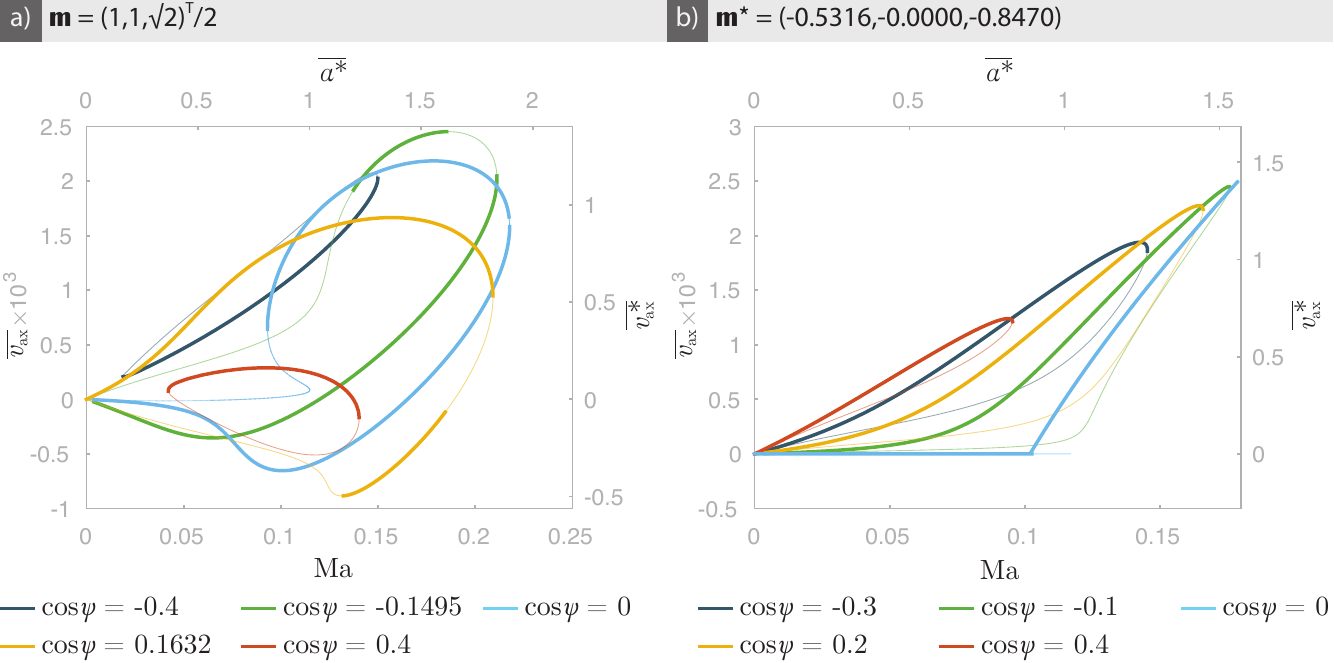}
	{\phantomsubcaption\label{subfig: morozov 3b m3 optangle}}
	{\phantomsubcaption\label{subfig: morozov 3b mstar optangle}}
	\caption[]{Three-bead cluster of \cite{Morozov2017} with optimal vertex angle, for several magnetic moment directions $\mbody$: axial velocities plotted against Mason number for different values of $\psi$. Thin lines correspond to unstable relative equilibria. Two scales are indicated: $\Ma$ and $\overline \vax$ refer to the Mason number and non-dimensional axial velocity as introduced in section~\ref{section: dim analysis}, while $\overline{\alpha^*}$ and $\overline{v_\text{ax}^*}$ refer to the non-dimensional angular velocity of the magnetic field and non-dimensional axial velocity using the scaling proposed by~\cite{Morozov2017}.}
	\label{fig: morozov 3b comp optangle}
\end{figure}
\begin{figure}[p!]
	\includegraphics{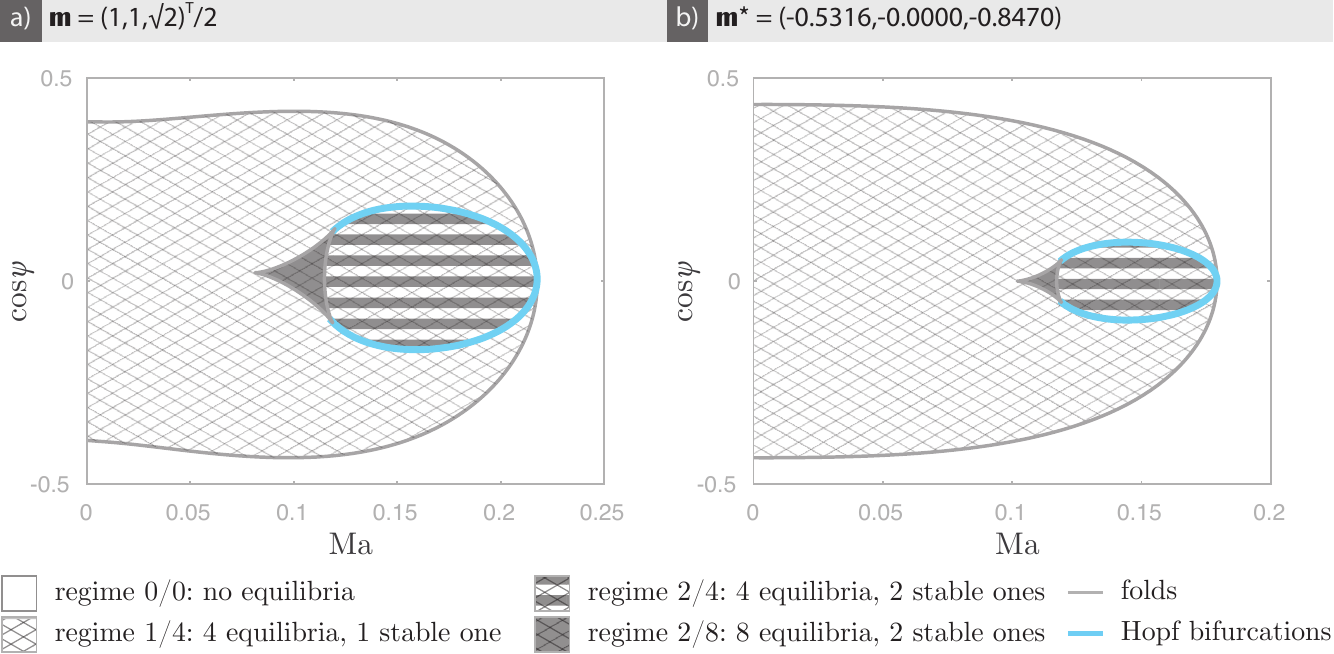}
	{\phantomsubcaption\label{subfig: morozov 3b m3 optangle regime}}
	{\phantomsubcaption\label{subfig: morozov 3b mstar optangle regime}}
	\caption[]{Three-bead cluster of \cite{Morozov2017} with optimal vertex angle, for several magnetic moment directions $\mbody$: diagram of the experimental parameters regimes.}
	\label{fig: morozov 3b regimes optangle}
\end{figure}

\begin{figure}[h!]
	\includegraphics{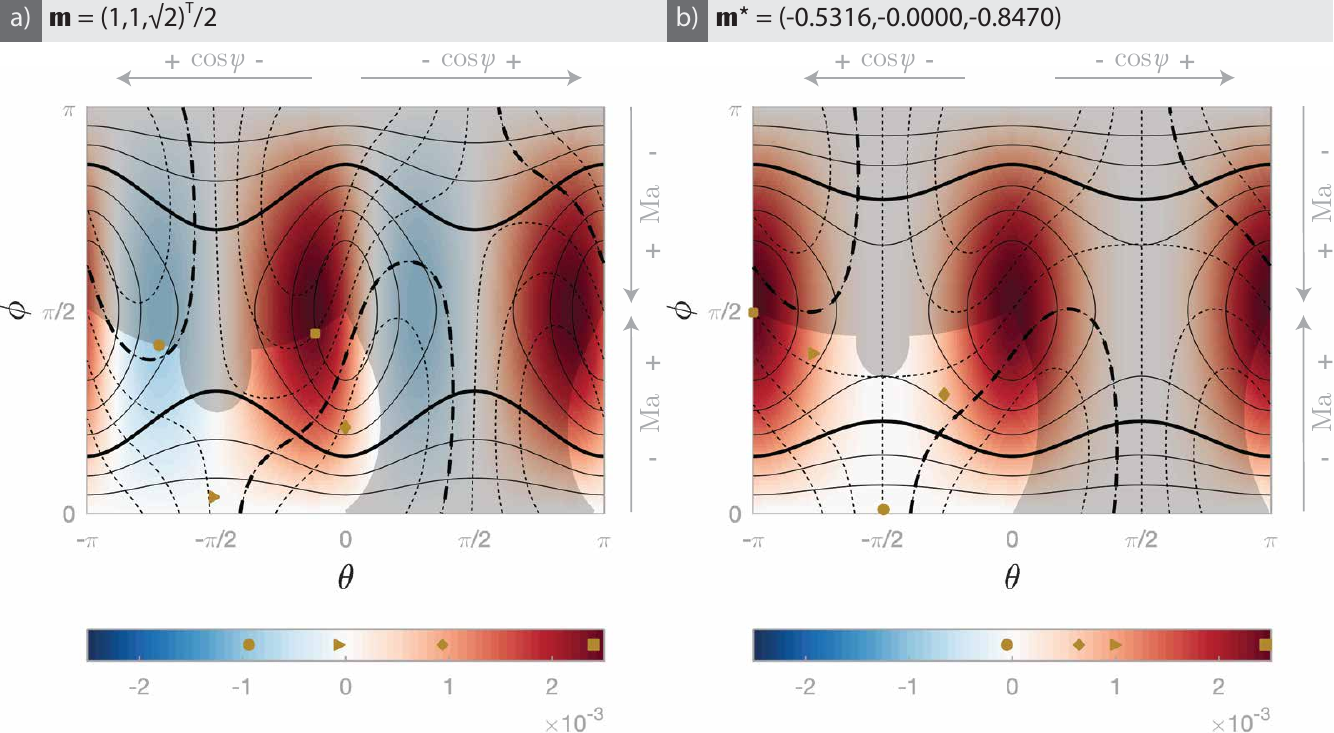}
	{\phantomsubcaption\label{subfig: morozov 3b m3 optangle contours}}
	{\phantomsubcaption\label{subfig: morozov 3b mstar optangle contours}}
	\caption[]{Three-bead cluster of \cite{Morozov2017} with optimal vertex angle, for several magnetic moment directions $\mbody$: axial velocity as a function of the mapping parameters $\theta$ and $\phi$ with level curves of the experimental parameters $\Ma$ in solid lines and $\cos\psi$ in dashed lines. Highlighted contours are \subref{subfig: morozov 3b m3 optangle contours}) $\Ma = 0.0933$ and $\cos \psi = 0.1338$, \subref{subfig: morozov 3b mstar optangle contours}) $\Ma = 0.0768$ and $\cos \psi = 0.1450$.}
	\label{fig: morozov 3b contours optangle}
\end{figure}

\cite{Morozov2017} also investigated how to optimise the axial velocity by varying the shape of swimmers. In particular, they studied three beads clusters with different vertex angles between the lines joining their centres, and found that the axial velocity is optimised for a vertex angle of $\simeq 122.7^\circ$. Using the mobility matrix they reported for this optimal swimmer (cf. appendix~\ref{appendix: numerics}), we show in figure~\ref{fig: morozov 3b comp optangle} the axial velocity as a function of the Mason number at different values of $\cos \psi$ for two distinct magnetisation directions: $\mbody = (1,1,\sqrt{2})^T/2$, and the optimal magnetisation direction $\mbody^\star = (-0.5316,-0.0000,-0.8470)^T$. Note that the largest axial velocity $\overline{\vax} = 0.0025$ is reached for $\cos \psi = 0$ with the optimal magnetisation direction, but an axial velocity only 2\% smaller is reached with $\cos \psi = -0.1495$ for the non optimal magnetisation direction $\mbody = (1,1,\sqrt{2})^T/2$. The corresponding experimental regimes diagram are provided in figure~\ref{fig: morozov 3b regimes optangle}, and figure~\ref{fig: morozov 3b contours optangle} shows the axial velocity and level sets of $\Ma$ and $\cos \psi$ as functions of the mapping parameters $\theta$ and $\phi$.

\section{Conclusions}

For all examples presented here, we did observe that when the magnetisation is optimal,  the optimal angle $\psi$ between the external field and its rotation axis was indeed $\pi/2$. In that case, we also observed a linear relation between the axial velocity and the Mason number (which linearly parameterises the angular velocity of the external field) over the range of Mason numbers on which the swimmer is bistable: the optimal equilibrium being reached just before step-out. Note that we do not know of a proof of this feature as a general result.

However, it is practically difficult to accurately control the magnetisation direction. For arbitrary magnetisation directions, optimal driving parameters were often found both well away from step-out and well away from $\psi=\pi/2$. In those cases, the relation between $\vax$ and $\Ma$ is strongly nonlinear. The good news however is that the optimal axial velocities achieved for all magnetisations were comparable in order of magnitude to those obtained with optimal magnetisation. The conclusion is therefore that in practice, we do not have to magnetise the swimmer optimally provided that we learn how to pinpoint optimal $\Ma$ and $\psi$. 

Furthermore, optimal driving parameters were always found in a bistable regime. It would be very interesting to have an experimental verification of the practicality of the methods proposed here for systematically reaching the faster equilibrium. 

Note that we have studied a non-dimensional version of the axial velocity. By plugging the dimensions back in, it is straightforward to show that the dimensional axial velocity is proportional to the angular velocity of the external field (and not to $\Ma$). Therefore, in practice, given a stable equilibrium, we can always multiply the axial velocity by a factor $\gamma$ by increasing both $\alpha$ and the magnitude of the external field the same factor $\gamma$. By doing so, we stay on the exact same non-dimensional equilibrium.

Finally, the mobility matrix of a given experimental swimmer and its magnetisation are not easy to measure in practice. We believe that an interesting line of enquiry would be to make use of the formal parameterisation of the equilibria proposed here  to reverse engineer the broad features of the experimental parameters diagram. One could for instance infer the outer boundary of the region of existence of stable equilibria by tracking $\Ma$ and $\cos\psi$ at step-out. It should then be possible to recover interesting informations about the  matrix~$\Pbody$ of the swimmer \textit{a posteriori}. Note for instance that $\sigma_1$ can be directly inferred from the maximal mason number for which there exists an equilibrium. The shape of the $\vax(\Ma)$ curve could also be used to infer properties of $\Pbody$. Once these are known, the methods of section~\ref{section: analysis} can be used to infer optimal driving parameters of this particular swimmer.

\section*{Acknowledgements}
We thank Oscar Gonzalez for providing the code used for computing the drag matrices, and John H. Maddocks for his support. This work was funded by the Swiss National Science Foundation (SNF) project number 200021-156403.

Declaration of interest. The authors report no conflict of interest.

\appendix

\begin{figure}[b!]
	\centering
	\includegraphics{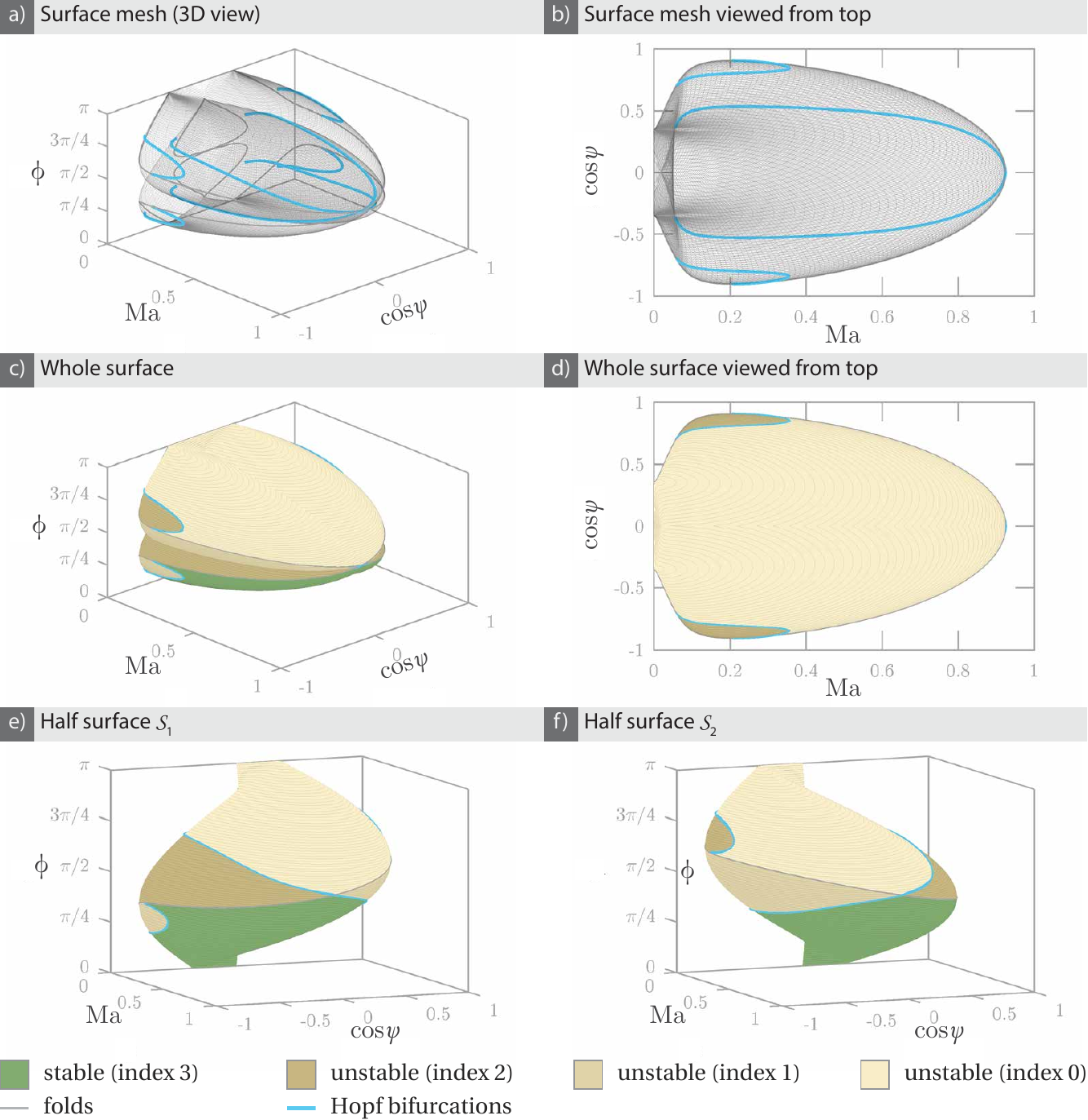}
	{\phantomsubcaption
	\label{subfig: swimmer B mesh}}
	{\phantomsubcaption
	\label{subfig: swimmer B mesh top}}
	{\phantomsubcaption
	\label{subfig: swimmer B surf}}
	{\phantomsubcaption
	\label{subfig: swimmer B surf top}}
	{\phantomsubcaption
	\label{subfig: swimmer B surf 1/2}}
	{\phantomsubcaption
	\label{subfig: swimmer B surf 2/2}}
	\caption[]{The analogue of fig.~\ref{fig: swimmer A eq} for swimmer B showing the surface~$\mathcal S$ corresponding to the set of equilibria of \eqref{ode Q closed} projected in the $\left( \Ma, \cos \psi, \phi \right)$-space.
	}
	\label{fig: swimmer B eq}
\end{figure}

\section{Numerical Data of Examples}
\label{appendix: numerics}

In section \ref{section: numerics}, swimmers A and B are helical rods with geometry and magnetisation directions as specified in table~\ref{table: swimmer geom}. The magnetisation directions are given in the basis in which the helical centreline of the rod is parametrised as
\begin{equation*}
	\begin{bmatrix} R \, \cos(s) \\ R \, \sin(s) \\ \frac{P}{2\pi} \end{bmatrix} \, ,
\end{equation*}
where $s = [0,s_\text{max}]$, $R$ is the helix radius, and $P$ is the helix pitch, so that the third basis direction corresponds to the helix axis. The total arc-length listed in table~\ref{table: swimmer geom} is the length of the helical centreline of the rod, and this rod is additionally capped by two half spheres with the same radius as the cross-section (cf.~fig.~\ref{fig: swimmers}).

\begin{table}[h]
	\centering
	\footnotesize{
	\begin{tabularx}{\textwidth}{l C C C C C}
		\toprule
		Helical swimmer & helix radius & helix pitch & total arc-length & cross-section radius & magnetisation direction ($\mbody$)\\
		\hline
		Swimmer A & $0.7158 \, \ell$ & $1.3789 \, \ell$ & $7.0566 \, \ell$ & $0.3343 \, \ell$ & $\begin{bsmallmatrix}\phantom{-}0.6755\\-0.7369\\\phantom{-}0.0245\end{bsmallmatrix}$ \\
		Swimmer B & $0.1330 \, \ell$ & $1.1076 \, \ell$ & $4.1628 \, \ell$ & $0.0936 \, \ell$ & $\begin{bsmallmatrix} 0\\0.1736\\0.9848 \end{bsmallmatrix}$\\
		\bottomrule
	\end{tabularx}}
	\caption{Specification of the geometry of helical swimmers.}
	\label{table: swimmer geom}
\end{table}

The drag matrices obtained from this information using Gonzalez' code~\citep[]{Gonzalez2009, Li2013, Gonzalez2015} are respectively
\begin{align*}
	&\dragbody_{11}^A = \eta \, \ell \, \begin{bsmallmatrix}
		20.8224  &  0.0000  &  0.5707 \\
   		-0.0000  & 21.0847  &  0.0000 \\
    		0.5708  &  0.0000  & 20.8330
	\end{bsmallmatrix} \, , &
	&\dragbody_{12}^A = \eta \, \ell^2 \, \begin{bsmallmatrix}
		0.5888  &  0.0000  &  0.5471 \\
   		0.0000 &  -0.0727  &  0.0000 \\
    		-0.0197 &   0.0000 &  -0.5366
	\end{bsmallmatrix} \, , \\
	&\dragbody_{22}^A = \eta \, \ell^3 \, \begin{bsmallmatrix}
		38.6928 &  -0.0000  &  4.6310 \\
		-0.0000  & 42.2129  &  0.0000 \\
		4.6312  &  0.0000  & 31.7645
	\end{bsmallmatrix} \, , &
	&\text{and} \\
	&\dragbody_{11}^B = \eta \, \ell \, \begin{bsmallmatrix}
		2.4654  &  0.0000  &  0.0000 \\
		-0.0000  & 12.4815  & 0.0582 \\
		-0.0000  &  0.0577   & 9.2808
	\end{bsmallmatrix} \, , &
	&\dragbody_{12}^B = \eta \, \ell^2 \, \begin{bsmallmatrix}
		0.1433  &  0.0000 &  -0.0000 \\
		0.0000  &  0.0122  &  0.1178 \\
		0.0000 &  -0.5607  & -0.2158
	\end{bsmallmatrix} \, , \\
	&\dragbody_{22}^B = \eta \, \ell^3 \, \begin{bsmallmatrix}
		20.1070  & -0.0000  &  0.0000 \\
		-0.0000  & 20.1725  &  0.4032 \\
		0.0000  &  0.4031  &  1.0196
	\end{bsmallmatrix} \, .
\end{align*}
Since these results are not exactly symmetric, we correct this by using $\Mob = \frac{1}{2} \Drag^{-1} + \frac{1}{2} \Drag^{-T}$ in our computations. Relevant material parameters are listed in table~\ref{table: swimmer params}.

\begin{table}[h]
	\centering
	\footnotesize{
	\begin{tabularx}{\textwidth}{l C C C C C C C}
		\toprule
		Swimmer & $\sigma_{1}$ & $\sigma_{2}$ & $c_{01}$ & $c_{02}$ & $c_{11}$ & $c_{12}$ & $\theta_{0}$ (rad)\\
		\hline
		Swimmer A & 0.0333 & 0.0241 & -0.0027 & $-3.6064$ $\times10^{-4}$ & $1.6878 $ $\times 10^{-5}$ & 0.0091 & -1.3871 \\
		Swimmer B & 0.9244 & 0.0497 & 0.3243 & $4.7998 $ $\times 10^{-5}$ & $-1.7003 $ $\times 10^{-5}$ & 0.8160 & -1.5680 \\
		\cite{Meshkati2014} & 0.1858 & 0.1274 & -0.0377 & 0.0095 & -0.0012 & 0.0549 & 1.2170 \\
		\cite{Morozov2017} \\
		\,\,Vertex angle: 90$^\circ$ \\
		$\mbody = (0,1,1)^T/\sqrt{2}$ & 0.1761 & 0.1190 & 0.0363 & $-2.3084 $ $\times 10^{-18}$ & $6.4068$ $\times10^{-19}$ & 0.0533 & 1.5708 \\
		$\mbody = (1,1,\sqrt{2})^T/2$ & 0.1735 & 0.1271 & 0.0397 & -0.0105 & -0.0014 & 0.0422 & 1.2253 \\
		$\mbody = (1,0,1)^T/\sqrt{2}$ & 0.1698 & 0.1360 & 0.0448 & $5.8280 $ $\times 10^{-20}$ & $-2.4619$ $\times 10^{-18}$ & -0.0278 & -1.5708\\
		$\mbody^\star = \begin{bsmallmatrix}-0.6026\\-0.0000\\-0.7980\end{bsmallmatrix}$ & 0.1575 & 0.1360 & -0.0431 & $7.4439 $ $\times 10^{-10}$ & $1.1162$ $\times 10^{-10}$ & -0.0156 & 1.5708 \\
		\,\,Vertex angle: 122.7$^\circ$ \\
		$\mbody = (1,1,\sqrt{2})^T/2$ & 0.2177 & 0.1148 & 0.0853 & -0.0026 & $-7.1608$ $\times 10^{-4}$ & 0.0854 & 1.5122 \\
		$\mbody^\star = \begin{bsmallmatrix}-0.5316\\-0.0000\\-0.8470\end{bsmallmatrix}$ & 0.1792 & 0.1172 & -0.0779 & $1.3648$ $\times 10^{-11}$ & $3.8181$ $\times 10^{-12}$ & -0.0442 & 1.5708 \\
		\bottomrule
	\end{tabularx}}
	\caption{Key material parameters for computing the relative equilibria}
	\label{table: swimmer params}
\end{table}

Three beads cluster with a vertex angle of 90$^\circ$ (taken from~\cite{Meshkati2014} and scaled taking $\ell = 2 \sqrt{2} \, R_{b} = 6.1518~\mu$m ($R_{b}$ is the radius of one bead) and assuming that the fluid viscosity is $\eta = 10^{-3}$~Pa s):
\begin{align*}
	&\mobbody_{11}^\text{3B 90$^\circ$} = \frac{1}{\eta \, \ell} \, \begin{bsmallmatrix}
		0.0849 & 0 & 0 \\
		0 & 0.1064 & 0 \\
		0 & 0 & 0.1058
	\end{bsmallmatrix} \, , &
	&\mobbody_{12}^\text{3B 90$^\circ$} = \frac{1}{\eta \, \ell^2} \, \begin{bsmallmatrix}
		0 & 0 & 0 \\
		0 & 0 & -0.0439 \\
		0 & 0.0734 & 0
	\end{bsmallmatrix} \, , \\
	&\mobbody_{22}^\text{3B 90$^\circ$} = \frac{1}{\eta \, \ell^3} \, \begin{bsmallmatrix}
		0.1360 & 0 & 0 \\
		0 & 0.2086 & 0 \\
		0 & 0 & 0.1190
	\end{bsmallmatrix} \, , &
	&\mbody = \frac{1}{\sqrt{3}} \, \begin{bsmallmatrix} 1 \\ 1 \\ 1 \end{bsmallmatrix} \, .
\end{align*}
This data is used to generate figures~\ref{fig: meshkati comp} and~\ref{fig: meshkati regimes+contours}.

Same matrices transformed according to the rules of \cite{Morozov2017}:
\begin{align*}
	&\mobbody_{12}^\text{3B 90$^\circ$ alt} = \frac{1}{\eta \, \ell^2} \, \begin{bsmallmatrix}
		0 & 0 & -0.0013 \\
		0 & 0 & 0 \\
		-0.0013 & 0 & 0
	\end{bsmallmatrix} \, , &
	&\mobbody_{22}^\text{3B 90$^\circ$ alt} = \frac{1}{\eta \, \ell^3} \, \begin{bsmallmatrix}
		 0.1190 & 0 & 0 \\
		0 & 0.1360 & 0 \\
		0 & 0 & 0.2086
	\end{bsmallmatrix} \, .
\end{align*}
Additional parameters used in the conversion between different scales (cf. fig.~\ref{fig: morozov 3b comp}): $|\mathrm{Ch}|  = |(\mobbody_{12})_{13} \, \big( (\mobbody_{22})_{11}^{-1} + (\mobbody_{22})_{33}^{-1} \big)/2| = 0.0085$ and $F_{\bot} = 2/\big( (\mobbody_{22})_{11}^{-1} + (\mobbody_{22})_{22}^{-1} \big) = 0.1269$. This data is used to generate figures~\ref{fig: morozov 3b comp}, \ref{fig: morozov 3b regimes}, and~\ref{fig: morozov 3b contours}.


Three beads cluster with optimal vertex angle of 122.7$^\circ$, taken from~\cite{Morozov2017}, and scaled taking $\ell = 2\sqrt{2} R_{b}$, where $R_{b}$ is the radius of one single bead:
\begin{align*}
	&\mobbody_{12}^\text{3B 122.7$^\circ$} = \frac{1}{\eta \, \ell^2} \, \begin{bsmallmatrix}
		0 & 0 & -0.0025 \\
		0 & 0 & 0 \\
		-0.0025 & 0 & 0
	\end{bsmallmatrix} \, , &
	&\mobbody_{22}^\text{3B 122.7$^\circ$}= \frac{1}{\eta \, \ell^3} \, \begin{bsmallmatrix}
		 0.1125 & 0 & 0 \\
		0 & 0.1172 & 0 \\
		0 & 0 & 0.2856
	\end{bsmallmatrix} \, .
\end{align*}
Additional parameters used in the conversion between different scales: $|\mathrm{Ch}|  = 0.0155$ and $F_{\bot} = 0.1148$. This data is used to generate figures~\ref{fig: morozov 3b comp optangle}, \ref{fig: morozov 3b regimes optangle}, and~\ref{fig: morozov 3b contours optangle}.


\section{More on the Geometry of the Surface of Equilibria}
\label{appendix: intersections}

\begin{figure}[b!]
	\centering
	\includegraphics{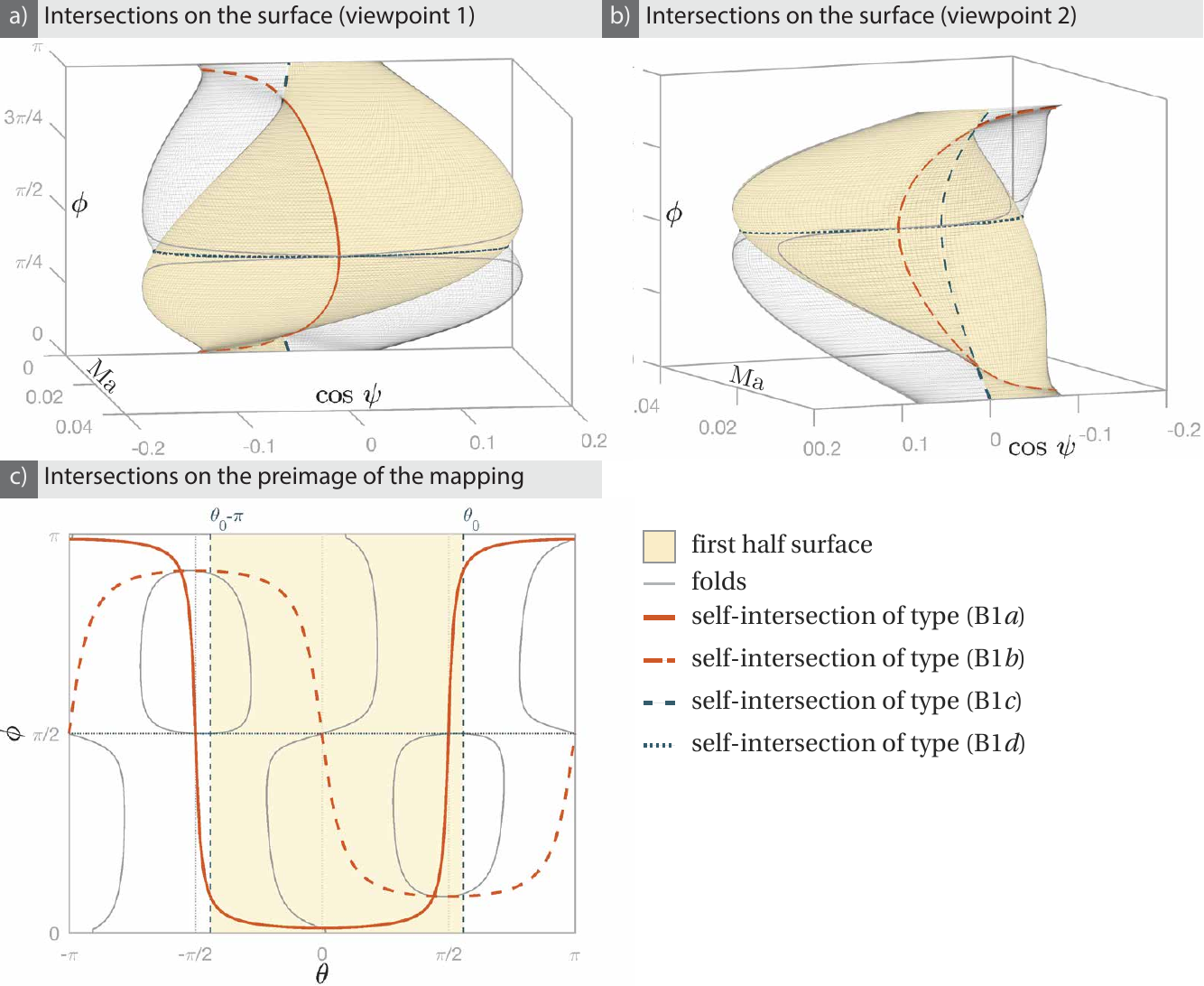}
	{\phantomsubcaption\label{subfig: swimmer A inters surf}}
	{\phantomsubcaption\label{subfig: swimmer A inters surf 2}}
	{\phantomsubcaption\label{subfig: swimmer A inters theta phi}}
	\caption[]{Self-intersections of the surface $\mathcal{S}$ shown \subref{subfig: swimmer A inters surf}-\subref{subfig: swimmer A inters surf 2}) on the surface and \subref{subfig: swimmer A inters theta phi}) on the plane of coordinates $\theta$ and $\phi$.}
	\label{fig: intersections}
\end{figure}

In order to better understand the geometry of the surface $\mathcal{S}$ described in section~\ref{sec-releq}, we give here a detailed account of all its self-intersections. The surface $\mathcal{S}$ is parametrised by variables $\theta$ and $\phi$. In section~\ref{section: surf charac} we indicated the self-intersection at $\theta = \theta_{0}$, where $\theta_{0}$ depends only on the swimmer, and used it to split the surface into two symmetric halves $\mathcal{S}_{1}$ and $\mathcal{S}_{2}$. Upon close inspection of the parametrisation~\eqref{param a psi} we find that the self-intersections of $\mathcal{S}$ are given by
\begin{subequations}
	\begin{align}
		\phi &= \mathrm{arccot} \left( - \frac{c_{12} \cos \theta}{c_{02} \, \sigma_{1}} \left( \frac{\cos^2 \theta}{\sigma_{1}^2} + \frac{\sin^2 \theta}{\sigma_{2}^2} \right)^{-1/2} \right) && \text{where } \Sigma \left( \theta, \phi \right) = \Sigma \left( - \theta, \phi \right) \label{inters theta=-theta}\\
		\phi &= \mathrm{arccot} \left( - \frac{c_{12} \sin \theta}{c_{01} \, \sigma_{2}} \left( \frac{\cos^2 \theta}{\sigma_{1}^2} + \frac{\sin^2 \theta}{\sigma_{2}^2} \right)^{-1/2} \right) && \text{where } \Sigma \left( \theta, \phi \right) = \Sigma \left( \pi - \theta, \phi \right) \label{inters theta=pi-theta} \\
		\theta &= \theta_{0}, \theta_{0} \pm \pi && \text{where } \Sigma \left( \theta, \phi \right) = \Sigma \left( \theta + \pi, \phi \right) \label{inters theta0} \\
		\phi &= \frac{\pi}{2} && \text{where } \Sigma \left( \theta, \phi \right) = \Sigma \left( \theta + \pi, \phi \right) \label{inters phi=pi/2}
	\end{align}
\end{subequations}

We can see by comparing the different panels of figures~\ref{fig: swimmer A eq} and~\ref{fig: swimmer B eq} that the two halves $\mathcal{S}_{1}$ and $\mathcal{S}_{2}$ intersect each other at multiple locations. A crucial question to understand the visualisation is whether there are self-intersections within the same half. We can readily see that there are on fig.~\ref{fig: swimmer A eq}. Intersection\re{inters theta0} defines the splitting between $\mathcal{S}_{1}$ and $\mathcal{S}_{2}$ while\re{inters phi=pi/2} is necessarily an intersection between the two halves by definition of $\mathcal{S}_{1}$ and $\mathcal{S}_{2}$. Therefore self-intersections within the same half are either in the family\re{inters theta=-theta} or\re{inters theta=pi-theta}.

\begin{figure}[b!]
	\centering
	\includegraphics{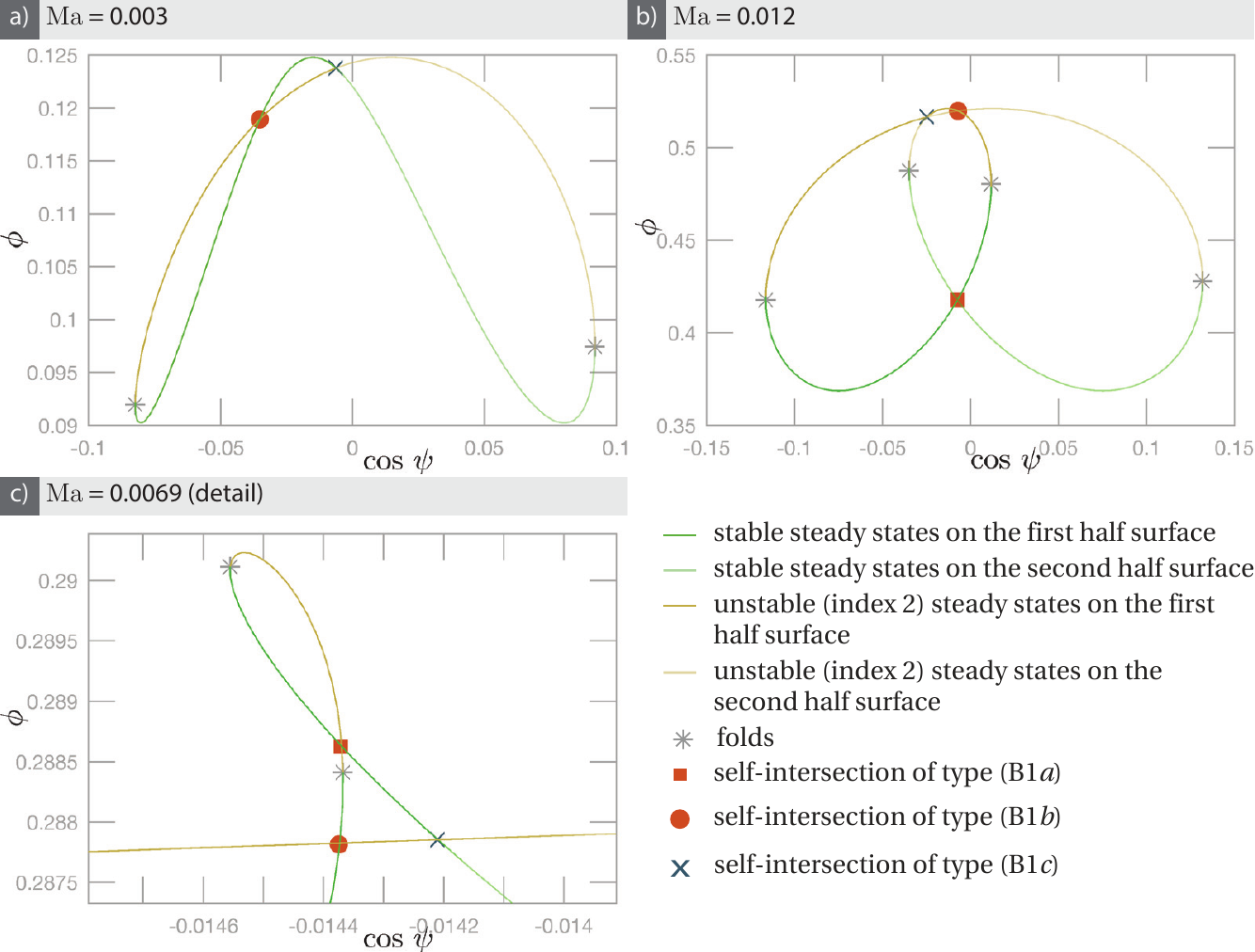}
	{\phantomsubcaption\label{subfig: a=0,003}}
	{\phantomsubcaption\label{subfig: a=0,012}}
	{\phantomsubcaption\label{subfig: a=0,0069}}
	\caption[]{Sections of surface $\mathcal{S}$ at constant $\Ma$. Panels \subref{subfig: a=0,003}-\subref{subfig: a=0,012}) show only section of the lower part of the surface ($\phi < \pi/2$); the upper part is symmetric to it. Panel \subref{subfig: a=0,0069} shows a detail of the lower part section, aimed at clarifying how section in panel~\subref{subfig: a=0,003} transforms into section in panel~\subref{subfig: a=0,012} as $\Ma$ grows. The data corresponds to swimmer A (cf. fig.~\ref{subfig: swimmer A inters surf})}
	\label{fig: a sections}
\end{figure}

Figure~\ref{fig: intersections} shows all the self-intersections of the surface $\mathcal S$ for swimmer A. On figure~\ref{subfig: swimmer A inters theta phi}, the type~\eqref{inters theta=-theta} intersection curve is symmetric about $\theta = 0$ and $\theta = \pm \pi$, while the type~\eqref{inters theta=pi-theta} intersection curve is symmetric about $\theta = \pm \frac{\pi}{2}$. Two points of either intersection curve in the $\left( \theta,\phi \right)$-plane coincide on $\mathcal{S}$ if they have the same height $\phi$. This results in self-intersections within the same half-surface $\mathcal{S}_{i}$ if those two points lie within the same region of the $\left( \theta, \phi \right)$-plane delimited by the $\theta = \theta_{0}$ and $\theta = \theta_{0} \pm \pi$ vertical lines (taking $\theta$-periodicity of the parametrisation into account), and in intersections between $\mathcal{S}_{1}$ and $\mathcal{S}_{2}$ otherwise.

Self-intersections of type \eqref{inters theta=-theta} of a half-surface $\mathcal{S}_{i}$ occur when $\theta$ is such that $\left| \tan \theta \right| < \left| \tan \theta_{0} \right|$: the closer $\theta_{0}$ is to $0$ or $\pm \pi$, the smaller the set of such $\theta$. Type \eqref{inters theta=pi-theta} self-intersections of $\mathcal{S}_{i}$ arise for $\theta$ such that $\left| \tan \theta \right| > \left| \tan \theta_{0} \right|$: the closer $\theta_{0}$ is to $\pm \frac{\pi}{2}$, the smaller the set of such $\theta$. For swimmer A, $\theta_{0}$ is quite close to $\frac{\pi}{2}$ (cf fig.~\ref{subfig: swimmer A inters theta phi}); therefore type~\eqref{inters theta=pi-theta} self-intersections within the same half $\mathcal{S}_{i}$ are indiscernible on figure~\ref{subfig: swimmer A inters surf}-\subref{subfig: swimmer A inters surf 2}). A larger set of $\theta$ satisfy the condition for type~\eqref{inters theta=-theta} self-intersections to occur within the same half $\mathcal{S}_{i}$.

Knowing how the surface self-intersects helps understanding how the half surfaces $\mathcal{S}_{1}$ and $\mathcal{S}_{2}$ join, and how the stability index is continuous across the two halves. We illustrate this by showing curves obtained as sections of $\mathcal{S}$ at constant values of $a$ for swimmer A in fig.~\ref{fig: a sections}.

\bibliographystyle{abbrvnat}
\bibliography{library}

\end{document}

%% file: notation.tex
\usepackage{isomath}
\usepackage{bm}
\usepackage{amssymb}
\SetMathAlphabet{\mathbf}{normal}{OML}{mdbch}{b}{n}

\def\dlab{\bm{d}}

\def\elab{\bm{e}}
\def\ebody{\bm{\mathsf{e}}}

\def\ulab{\bm{u}}
\def\ubody{\bm{\mathsf{u}}}

\def\xlab{\bm{x}}

\def\vlab{\bm{v}}
\def\vbody{\bm{\mathsf{v}}}
\def\vbodybar{\overline{\vbody}}
\def\omegalab{\bm{\omega}}
\def\omegabody{\mathbf{\omega}} 
\def\omegabodybar{\overline{\omegabody}}

\def\inertlab{I_{\mathrm{cm}}}
\def\inertbody{\mathsf{I}_{\mathrm{cm}}}
\def\inertbodybar{\overline{\inertbody}}

\def\plab{\bm{p}}
\def\pbody{\bm{\mathsf{p}}}
\def\pbodybar{\overline{\pbody}}

\def\Llab{\bm{L}}
\def\Lbody{\bm{\mathsf{L}}}
\def\Lbodybar{\overline{\Lbody}}

\def\flab{\bm{f}}
\def\fbody{\bm{\mathsf{f}}}
\def\fbodybar{\overline{\fbody}}
\def\taulab{\bm{\tau}}
\def\taubody{\mathbf{\tau}}
\def\taubodybar{\overline{\taubody}}
\def\msup{^{\left( m \right)}}
\def\dsup{^{\left( d \right)}}
\def\bsup{^{\left( b \right)}}

\def\extsup{^{\left( \text{ext} \right)}}

\def\Drag{\mathbb{D}}
\def\draglab{D}
\def\dragbody{\mathsf{D}}
\def\dragbodybar{\overline{\dragbody}}
\def\Mob{\mathbb{M}}

\def\mobbody{\mathsf{M}}

\def\gbody{\bm{\mathsf{g}}}

\def\gbodybar{\overline{\gbody}}

\def\Deltabody{\mathbf{{\Delta}}}
\def\Deltabodybar{\overline{\Deltabody}}

\def\mlab{\bm{m}}
\def\mbody{\bm{\mathsf{m}}}
\def\mbodybar{\overline{\mbody}}

\def\Blab{\bm{B}}
\def\Bbody{\bm{\mathsf{B}}}
\def\Bbodybar{\overline{\Bbody}}

\def\nbody{\boldsymbol{\mathsf{n}}}

\def\ubody{\bm{\mathsf{u}}}

\def\tbar{\overline{t}}

\def\Pbody{\mathsf{P}}
\def\betabody{\mathbf{\beta}}
\def\etabody{\mathbf{\eta}}

\def\Rey{\mathrm{Re}}
\def\eps{\varepsilon}

\usepackage{ltxcmds}
\makeatletter
\newcommand{\re}[1]{~(\ref{#1}\checknextarg}
\newcommand{\checknextarg}{\ltx@ifnextchar@nospace\bgroup{\gobblenextarg}{)}}
\newcommand{\gobblenextarg}[1]{,~\ref{#1}\ltx@ifnextchar@nospace\bgroup{\gobblenextarg}{)}}
\makeatother

\renewcommand{\re}[1]{~\eqref{#1}}
\newcommand{\ret}[2]{~(\ref{#1},~\ref{#2})}

\newcommand{\bma}{\begin{bmatrix}}
\newcommand{\ema}{\end{bmatrix}}

\def\Ma{\mathrm{Ma}}

\def\vax{v_\text{ax}}

%% file: draft.bbl
\begin{thebibliography}{38}
\providecommand{\natexlab}[1]{#1}
\providecommand{\url}[1]{\texttt{#1}}
\expandafter\ifx\csname urlstyle\endcsname\relax
  \providecommand{\doi}[1]{doi: #1}\else
  \providecommand{\doi}{doi: \begingroup \urlstyle{rm}\Url}\fi

\bibitem[Bell et~al.(2007)Bell, Leutenegger, Hammar, Dong, and
  Nelson]{Bell2007}
D.~J. Bell, S.~Leutenegger, K.~M. Hammar, L.~X. Dong, and B.~J. Nelson.
\newblock Flagella-like {{Propulsion}} for {{Microrobots Using}} a {{Nanocoil}}
  and a {{Rotating Electromagnetic Field}}.
\newblock In \emph{Proceedings 2007 {{IEEE International Conference}} on
  {{Robotics}} and {{Automation}}}, pages 1128--1133, Apr. 2007.
\newblock \doi{10.1109/ROBOT.2007.363136}.

\bibitem[Bente et~al.(2018)Bente, Codutti, Bachmann, and Faivre]{Bente2018}
K.~Bente, A.~Codutti, F.~Bachmann, and D.~Faivre.
\newblock Biohybrid and {{Bioinspired Magnetic Microswimmers}}.
\newblock \emph{Small}, 14\penalty0 (29):\penalty0 1704374, July 2018.
\newblock ISSN 1613-6829.
\newblock \doi{10.1002/smll.201704374}.

\bibitem[Crenshaw(1993)]{Crenshaw1993}
H.~C. Crenshaw.
\newblock Orientation by {{Helical Motion}}--{{I}}. {{Kinematics}} of the
  {{Helical Motion}} of {{Organisms}} with up to six {{Degrees}} of
  {{Freedom}}.
\newblock \emph{Bulletin of Mathematical Biology}, 55\penalty0 (1):\penalty0
  197--212, 1993.

\bibitem[Dhooge et~al.(2008)Dhooge, Govaerts, Kuznetsov, Meijer, and
  Sautois]{Dhooge2008}
A.~Dhooge, W.~Govaerts, Y.~A. Kuznetsov, H.~G.~E. Meijer, and B.~Sautois.
\newblock New features of the software {{MatCont}} for bifurcation analysis of
  dynamical systems.
\newblock \emph{Mathematical and Computer Modelling of Dynamical Systems},
  14\penalty0 (2):\penalty0 147--175, Apr. 2008.
\newblock ISSN 1387-3954.
\newblock \doi{10.1080/13873950701742754}.

\bibitem[Dichmann et~al.(1996)Dichmann, Li, and Maddocks]{Dichmann1996}
D.~J. Dichmann, Y.~Li, and J.~H. Maddocks.
\newblock Hamiltonian {{Formulations}} and {{Symmetries}} in {{Rod Mechanics}}.
\newblock \emph{Mathematical Approaches to Biomolecular Structure and Dynamics,
  IMA Volumes in Mathematics and Its Applications}, 82:\penalty0 71--113, 1996.

\bibitem[Dreyfus et~al.(2005)Dreyfus, Baudry, Roper, Fermigier, Stone, and
  Bibette]{Dreyfus2005}
R.~Dreyfus, J.~Baudry, M.~L. Roper, M.~Fermigier, H.~A. Stone, and J.~Bibette.
\newblock Microscopic artificial swimmers.
\newblock \emph{Nature}, 437\penalty0 (7060):\penalty0 862--865, Oct. 2005.
\newblock ISSN 1476-4687.
\newblock \doi{10.1038/nature04090}.

\bibitem[Ebbens and Howse(2010)]{Ebbens2010}
S.~J. Ebbens and J.~R. Howse.
\newblock In pursuit of propulsion at the nanoscale.
\newblock \emph{Soft Matter}, 6\penalty0 (4):\penalty0 726--738, 2010.
\newblock \doi{10.1039/B918598D}.

\bibitem[Felfoul et~al.(2016)Felfoul, Mohammadi, Taherkhani, {de Lanauze}, Xu,
  Loghin, Essa, Jancik, Houle, Lafleur, Gaboury, Tabrizian, Kaou, Atkin, Vuong,
  Batist, Beauchemin, Radzioch, and Martel]{Felfoul2016}
O.~Felfoul, M.~Mohammadi, S.~Taherkhani, D.~{de Lanauze}, Y.~Z. Xu, D.~Loghin,
  S.~Essa, S.~Jancik, D.~Houle, M.~Lafleur, L.~Gaboury, M.~Tabrizian, N.~Kaou,
  M.~Atkin, T.~Vuong, G.~Batist, N.~Beauchemin, D.~Radzioch, and S.~Martel.
\newblock Magneto-aerotactic bacteria deliver drug-containing nanoliposomes to
  tumour hypoxic regions.
\newblock \emph{Nature Nanotechnology}, 11\penalty0 (11):\penalty0 941--947,
  Nov. 2016.
\newblock ISSN 1748-3395.
\newblock \doi{10.1038/nnano.2016.137}.

\bibitem[Fu et~al.(2015)Fu, Jabbarzadeh, and Meshkati]{Fu2015}
H.~C. Fu, M.~Jabbarzadeh, and F.~Meshkati.
\newblock Magnetization directions and geometries of helical microswimmers for
  linear velocity-frequency response.
\newblock \emph{Physical Review E}, 91\penalty0 (4):\penalty0 1--13, 2015.
\newblock \doi{10.1103/PhysRevE.91.043011}.

\bibitem[Gao and Wang(2014)]{Gao2014}
W.~Gao and J.~Wang.
\newblock The {{Environmental Impact}} of {{Micro}}/{{Nanomachines}}: {{A
  Review}}.
\newblock \emph{ACS Nano}, 8\penalty0 (4):\penalty0 3170--3180, Apr. 2014.
\newblock ISSN 1936-0851.
\newblock \doi{10.1021/nn500077a}.

\bibitem[Ghosh and Fischer(2009)]{Ghosh2009}
A.~Ghosh and P.~Fischer.
\newblock Controlled {{Propulsion}} of {{Artificial Magnetic Nanostructured
  Propellers}}.
\newblock \emph{Nano Letters}, 9\penalty0 (6):\penalty0 2243--5, June 2009.
\newblock \doi{10.1021/nl900186w}.

\bibitem[Ghosh et~al.(2012)Ghosh, Paria, Singh, Venugopalan, and
  Ghosh]{Ghosh2012}
A.~Ghosh, D.~Paria, H.~J. Singh, P.~L. Venugopalan, and A.~Ghosh.
\newblock Dynamical configurations and bistability of helical nanostructures
  under external torque.
\newblock \emph{Physical Review E}, 86\penalty0 (3), Sept. 2012.
\newblock \doi{10.1103/PhysRevE.86.031401}.

\bibitem[Ghosh et~al.(2013)Ghosh, Mandal, Karmakar, and Ghosh]{Ghosh2013}
A.~Ghosh, P.~Mandal, S.~Karmakar, and A.~Ghosh.
\newblock Analytical theory and stability analysis of an elongated nanoscale
  object under external torque.
\newblock \emph{Physical Chemistry Chemical Physics}, 15\penalty0
  (26):\penalty0 10817--10823, 2013.
\newblock \doi{10.1039/C3CP50701G}.

\bibitem[Gonzalez(2009)]{Gonzalez2009}
O.~Gonzalez.
\newblock On stable, complete, and singularity-free boundary integral
  formulations of exterior {{Stokes}} flow.
\newblock \emph{SIAM Journal on Applied Mathematics}, 69\penalty0 (4):\penalty0
  933--958, 2009.

\bibitem[Gonzalez and Li(2015)]{Gonzalez2015}
O.~Gonzalez and J.~Li.
\newblock A convergence theorem for a class of {{Nystr{\"o}m}} methods for
  weakly singular integral equations on surfaces in {{R3}}.
\newblock \emph{Mathematics of Computation}, 84\penalty0 (292):\penalty0
  675--714, 2015.

\bibitem[Gonzalez et~al.(2004)Gonzalez, Graf, and Maddocks]{Gonzalez2004}
O.~Gonzalez, A.~B.~A. Graf, and J.~H. Maddocks.
\newblock Dynamics of a {{Rigid Body}} in a {{Stokes Fluid}}.
\newblock \emph{Journal of Fluid Mechanics}, 519:\penalty0 133--160, 2004.

\bibitem[Happel and Brenner(1983)]{Happel1983}
J.~Happel and H.~Brenner.
\newblock \emph{Low {{Reynolds}} Number Hydrodynamics: With Special
  Applications to Particulate Media}.
\newblock Mechanics of {{Fluids}} and {{Transport Processes}}. {Springer
  Netherlands}, 1983.
\newblock ISBN 978-90-247-2877-0.

\bibitem[Hinch(1991)]{Hinch1991}
E.~J. Hinch.
\newblock \emph{Perturbation {{Methods}}}.
\newblock {Cambridge University Press}, 1991.

\bibitem[Honda et~al.(1996)Honda, Arai, and Ishiyama]{Honda1996}
T.~Honda, K.~I. Arai, and K.~Ishiyama.
\newblock Micro {{Swimming Mechanisms Propelled}} by {{External Magnetic
  Fields}}.
\newblock \emph{IEEE Transactions on Magnetics}, 32\penalty0 (5):\penalty0
  5085--5087, 1996.

\bibitem[Kim and Karrila(2013)]{Kim2013}
S.~Kim and S.~J. Karrila.
\newblock \emph{Microhydrodynamics: Principles and Selected Applications}.
\newblock {Courier Corporation}, 2013.

\bibitem[Kuznetsov(2004)]{Kuznetsov2004}
Y.~A. Kuznetsov.
\newblock \emph{Elements of Applied Bifurcation Theory}.
\newblock Number 112 in Applied Mathematical Sciences. {Springer}, {New York,
  NY}, 3. ed edition, 2004.
\newblock ISBN 978-1-4419-1951-9.
\newblock OCLC: 815949776.

\bibitem[Li and Gonzalez(2013)]{Li2013}
J.~Li and O.~Gonzalez.
\newblock Convergence and conditioning of a {{Nystr{\"o}m}} method for
  {{Stokes}} flow in exterior three-dimensional domains.
\newblock \emph{Advances in Computational Mathematics}, 39\penalty0
  (1):\penalty0 143--174, 2013.
\newblock \doi{10.1007/s10444-012-9272-1}.

\bibitem[Mahoney et~al.(2014)Mahoney, Nelson, Peyer, Nelson, and
  Abbott]{Mahoney2014}
A.~W. Mahoney, N.~D. Nelson, K.~E. Peyer, B.~J. Nelson, and J.~J. Abbott.
\newblock Behavior of rotating magnetic microrobots above the step-out
  frequency with application to control of multi-microrobot systems.
\newblock \emph{Applied Physics Letters}, 104\penalty0 (14):\penalty0 1--5,
  2014.
\newblock \doi{10.1063/1.4870768}.

\bibitem[Man and Lauga(2013)]{Man2013}
Y.~Man and E.~Lauga.
\newblock The wobbling-to-swimming transition of rotated helices.
\newblock \emph{Physics of Fluids}, 25\penalty0 (7):\penalty0 071904, July
  2013.
\newblock ISSN 1070-6631.
\newblock \doi{10.1063/1.4812637}.

\bibitem[{Medina-S{\'a}nchez} et~al.(2016){Medina-S{\'a}nchez}, Schwarz, Meyer,
  Hebenstreit, and Schmidt]{Medina-Sanchez2016}
M.~{Medina-S{\'a}nchez}, L.~Schwarz, A.~K. Meyer, F.~Hebenstreit, and O.~G.
  Schmidt.
\newblock Cellular {{Cargo Delivery}}: {{Toward Assisted Fertilization}} by
  {{Sperm}}-{{Carrying Micromotors}}.
\newblock \emph{Nano Letters}, 16\penalty0 (1):\penalty0 555--561, Jan. 2016.
\newblock ISSN 1530-6984.
\newblock \doi{10.1021/acs.nanolett.5b04221}.

\bibitem[Meshkati and Fu(2014)]{Meshkati2014}
F.~Meshkati and H.~C. Fu.
\newblock Modeling rigid magnetically rotated microswimmers: {{Rotation}} axes,
  bistability, and controllability.
\newblock \emph{Physical Review E}, 90\penalty0 (6):\penalty0 1--11, 2014.
\newblock \doi{10.1103/PhysRevE.90.063006}.

\bibitem[Meshkati and Fu(2017)]{Meshkati2017}
F.~Meshkati and H.~C. Fu.
\newblock Erratum: {{Modeling}} rigid magnetically rotated microswimmers:
  {{Rotation}} axes, bistability, and controllability [{{Phys}}. {{Rev}}. {{E}}
  {\textbf{90}} , 063006 (2014)].
\newblock \emph{Physical Review E}, 95\penalty0 (6), June 2017.
\newblock ISSN 2470-0045, 2470-0053.
\newblock \doi{10.1103/PhysRevE.95.069904}.

\bibitem[Morozov and Leshansky(2014)]{Morozov2014}
K.~I. Morozov and A.~M. Leshansky.
\newblock The chiral magnetic nanomotors.
\newblock \emph{Nanoscale}, 6\penalty0 (3):\penalty0 1580--1588, 2014.
\newblock \doi{10.1039/C3NR04853E}.

\bibitem[Morozov et~al.(2017)Morozov, Mirzae, Kenneth, and
  Leshansky]{Morozov2017}
K.~I. Morozov, Y.~Mirzae, O.~Kenneth, and A.~M. Leshansky.
\newblock Dynamics of arbitrary shaped propellers driven by a rotating magnetic
  field.
\newblock \emph{Physical Review Fluids}, 2\penalty0 (4):\penalty0 044202, Apr.
  2017.
\newblock \doi{10.1103/PhysRevFluids.2.044202}.

\bibitem[Mushtaq et~al.(2015)Mushtaq, Guerrero, Sakar, Hoop, M.~Lindo, Sort,
  Chen, J.~Nelson, Pellicer, and Pan{\'e}]{Mushtaq2015}
F.~Mushtaq, M.~Guerrero, M.~S. Sakar, M.~Hoop, A.~M.~Lindo, J.~Sort, X.~Chen,
  B.~J.~Nelson, E.~Pellicer, and S.~Pan{\'e}.
\newblock Magnetically driven {{Bi}} 2 {{O}} 3 /{{BiOCl}}-based hybrid
  microrobots for photocatalytic water remediation.
\newblock \emph{Journal of Materials Chemistry A}, 3\penalty0 (47):\penalty0
  23670--23676, 2015.
\newblock \doi{10.1039/C5TA05825B}.

\bibitem[Nelson et~al.(2010)Nelson, Kaliakatsos, and Abbott]{Nelson2010}
B.~J. Nelson, I.~K. Kaliakatsos, and J.~J. Abbott.
\newblock Microrobots for {{Minimally Invasive Medicine}}.
\newblock \emph{Annual Review of Biomedical Engineering}, 12:\penalty0 55--85,
  Aug. 2010.
\newblock \doi{10.1146/annurev-bioeng-010510-103409}.

\bibitem[Peyer et~al.(2013)Peyer, Zhang, and Nelson]{Peyer2013b}
K.~E. Peyer, L.~Zhang, and B.~J. Nelson.
\newblock Bio-{{Inspired Magnetic Swimming Microrobots}} for {{Biomedical
  Applications}}.
\newblock \emph{Nanoscale}, 5\penalty0 (4):\penalty0 1259--72, Feb. 2013.
\newblock \doi{10.1039/c2nr32554c}.

\bibitem[Powell(1981)]{Powell1981}
M.~J.~D. Powell.
\newblock \emph{Approximation Theory and Methods}.
\newblock {Cambridge University Press}, 1981.

\bibitem[{R{\"u}egg-Reymond}(2019)]{Ruegg-Reymond2019}
P.~{R{\"u}egg-Reymond}.
\newblock \emph{On the {{Dynamics}} of {{Magnetic Micro}}-{{Swimmers}}}.
\newblock PhD thesis, EPFL, 2019.

\bibitem[{R{\"u}egg-Reymond} and Lessinnes(2018)]{Ruegg-Reymond2018}
P.~{R{\"u}egg-Reymond} and T.~Lessinnes.
\newblock Asymptotic {{Dynamics}} of {{Magnetic Micro}}-{{Swimmers}}.
\newblock \emph{arXiv:1807.09059}, 2018.

\bibitem[Tottori et~al.(2012)Tottori, Zhang, Qiu, Krawczyk,
  {Franco-Obreg{\'o}n}, and Nelson]{Tottori2012}
S.~Tottori, L.~Zhang, F.~Qiu, K.~K. Krawczyk, A.~{Franco-Obreg{\'o}n}, and
  B.~J. Nelson.
\newblock Magnetic helical micromachines: Fabrication, controlled swimming, and
  cargo transport.
\newblock \emph{Advanced Materials}, 24\penalty0 (6):\penalty0 811--816, Feb.
  2012.
\newblock \doi{10.1002/adma.201103818}.

\bibitem[Vach et~al.(2015)Vach, Fratzl, Klumpp, and Faivre]{Vach2015}
P.~J. Vach, P.~Fratzl, S.~Klumpp, and D.~Faivre.
\newblock Fast {{Magnetic Micropropellers}} with {{Random Shapes}}.
\newblock \emph{Nano Letters}, 15\penalty0 (10):\penalty0 7064--7070, Oct.
  2015.
\newblock ISSN 1530-6984.
\newblock \doi{10.1021/acs.nanolett.5b03131}.

\bibitem[Yan et~al.(2017)Yan, Zhou, Vincent, Deng, Yu, Xu, Xu, Tang, Bian,
  Wang, Kostarelos, and Zhang]{Yan2017}
X.~Yan, Q.~Zhou, M.~Vincent, Y.~Deng, J.~Yu, J.~Xu, T.~Xu, T.~Tang, L.~Bian,
  Y.-X.~J. Wang, K.~Kostarelos, and L.~Zhang.
\newblock Multifunctional biohybrid magnetite microrobots for imaging-guided
  therapy.
\newblock \emph{Science Robotics}, 2\penalty0 (12):\penalty0 eaaq1155, Nov.
  2017.
\newblock ISSN 2470-9476.
\newblock \doi{10.1126/scirobotics.aaq1155}.

\end{thebibliography}
